\documentclass{amsart}
\usepackage{amsfonts}

\usepackage{amsthm}
\usepackage{amsfonts}
\usepackage{mathrsfs}
\usepackage{amscd,amssymb,amsmath,graphicx,verbatim,xypic,bm}
\usepackage{hyperref}
\usepackage[TS1,OT1,T1]{fontenc}

\newtheorem{theorem}{Theorem}[section]
\newtheorem{lemma}[theorem]{Lemma}

\newtheorem{proposition}[theorem]{Proposition}

\theoremstyle{definition}
\newtheorem{definition}[theorem]{Definition}

\theoremstyle{remark}
\newtheorem{remark}[theorem]{Remark}

\numberwithin{equation}{section}
\catcode`@=11
\let\@int\int \def\int{\displaystyle\@int}
\let\@lim\lim \def\lim{\displaystyle\@lim}
\let\@sum\sum \def\sum{\displaystyle\@sum}
\let\@sup\sup \def\sup{\displaystyle\@sup}
\let\@inf\inf \def\inf{\displaystyle\@inf}
\let\@cap\cap \def\cap{\displaystyle\@cap}
\let\@cup\cup \def\cup{\displaystyle\@cup}
\let\@max\max \def\max{\displaystyle\@max}
\let\@min\min \def\min{\displaystyle\@min}
\let\@frac\frac \def\frac{\displaystyle\@frac}
\let\@iint\iint \def\iint{\displaystyle\@iint}

\def\epsilon{\varepsilon}

\begin{document}

\title[Topologically Conjugate Classifications]
{Topologically Conjugate Classifications of the Translation Actions
on Compact Connected Lie Groups $\mathbf{SU(2)}\bm{\times}\bm{T^n}$}
\author{Xiaotian Pan}
\address{Xiaotian Pan, School of Mathematical Science, Capital Normal University, 100048, Beijing, P.R.China} \email{616280954@qq.com}

\author{Bingzhe Hou }
\address{Bingzhe Hou, Department of Mathematics, Jilin university, 130012, Changchun, P.R.China} \email{houbz@jlu.edu.cn}

\date{Mar. 6, 2022}
\begin{abstract}
In this article, we focus on the left (translation) actions on noncommutative compact connected Lie groups ${\rm SU}(2) \times T^n$.
We define the rotation vectors of the left actions induced by the elements in the maximal tori
of ${\rm SU}(2) \times T^n$,
and utilize rotation vectors to give  the complete topologically conjugate
classifications of left actions. Algebraic conjugacy and smooth conjugacy are also considered.
\end{abstract}
\maketitle

\section{Introduction}
The main research direction of this paper derives from S. Smale \cite{SS-1963}.
He hoped to classify all smooth self-maps of differential manifolds, and to obtain a
better understanding of the dynamical properties of some smooth self-maps.

In the present paper, we consider the topological case of this problem.
Assume that $X$ is a topological manifold, and $f,\,g$ are continuous self-maps of $X$.
If there exists a homeomorphism $h:\,X \rightarrow X$ such that
$$
h \circ f=g \circ h,
$$
$f$ and $g$ are said to be topologically conjugate.
Similarly, we have the concepts of smooth conjugacy and algebraic conjugacy.

Topologically conjugate classification is an important and difficult task in dynamical systems.
In the introduction of \cite{{HP-2018}}, we introduced some important works and results about topologically conjugate classifications, such as Anosov diffeomorphisms on differential manifolds \cite{JF-1970,MA-1974},
linear endomorphisms of finite dimensional vector spaces \cite{ATW-1997, KR-1973, RJW-1972},
linear representations of compact Lie groups \cite{SR-1977},
automorphisms of abelian topological groups \cite{AP-1965} and affine transformations of tori and more general compact abelian groups \cite{WP-1968,WP-1969} for non-minimal dynamical systems.
We want to extend the objects of this study,
and focus on the translation actions on some noncommutative compact connected Lie groups.
One can see that the translation actions on noncommutative compact connected Lie groups are all self-homeomorphisms, but they are neither affine transformations nor endomorphisms,
so our investigations are not contained in the dynamical systems mentioned above.
For translation actions on compact Lie groups,
S. Weinberger \cite{SW-2008} stated that $(G, g)$ is topologically conjugate to $(G, h)$ if and only if $G/{\rm cl}(\langle g \rangle) \cong G/{\rm cl}(\langle h \rangle)$ by an
isomorphism that pulls back the principle ${\rm cl}\langle g \rangle={\rm cl}\langle h \rangle$ bundles, where $<g>$ is the subgroup generated by $g$.
Our purpose is to find an explicit and computable way to classify the translation actions on a class of noncommutative compact connected Lie groups completely.
Following some investigations, we \cite{{HP-2018}} successfully utilized rotation vectors (numbers) to give the complete topologically conjugate classifications of the left actions on noncommutative compact connected Lie groups with topological dimension 3 or 4 .
So we hope to use the same way to study the left actions on the Lie groups ${\rm SU}(2) \times T^n$ in the present paper.

\begin{remark}
For convenience, we only investigate the left actions on the Lie groups ${\rm SU}(2) \times T^n$,
and the relevant conclusions of the right actions on these Lie groups are the same as left actions.
\end{remark}

Firstly, let us review the definitions of the rotations of $T^n$ and their rotation vectors.
For any
$$
\alpha=\left(\begin{array}{c}
{\textrm{e}}^{2\pi{\rm i}\theta_1}\\
{\textrm{e}}^{2\pi{\rm i}\theta_2}\\
\vdots\\
{\textrm{e}}^{2\pi{\rm i}\theta_n}
\end{array}\right) \in T^n,\qquad
\theta_1,\,\theta_2,\, \cdots \,{\theta}_{n}\in [0,1),
$$
we define a rotation $f_\alpha$ of $T^n$ by
$$
f_\alpha:\,u \rightarrow \left(\begin{array}{cccc}
{\textrm{e}}^{2\pi{\rm i}\theta_1}&&&\\
&{\textrm{e}}^{2\pi{\rm i}\theta_2}&&\\
&&\ddots&\\
&&&{\textrm{e}}^{2\pi{\rm i}\theta_n}
\end{array}\right)u,\qquad
\forall\, u\in T^n.
$$
Set
$$
\mathscr{M}_{T^n}=\big\{f_{\alpha}:\,T^n \rightarrow T^n\,|\,\alpha \in {T^n} \big\}.
$$
Then $\mathscr{M}_{T^n}$ is the set consisting of all rotations of $T^n$.

\begin{remark}
In the present paper, when we regard $T$ as a quotient space of $\mathbb{R}$ and talk about lift,
we always take the quotient map $\pi$ defined by $\pi(x)=\textrm{e}^{2\pi{\rm i}x},\,\forall x \in \mathbb{R}$.
\end{remark}

\begin{definition}\label{def:0}
Assume that $f_{\alpha} \in \mathscr{M}_{T^n}$, and $F_{\alpha}:\,\mathbb{R}^n \rightarrow \mathbb{R}^n$ is a lift of $f_{\alpha}$. Then the limit
$$
\rho(F_{\alpha})={\lim\limits_{n\to\infty}}{\dfrac{F_{\alpha}^{n}({\bm x})-{\bm x}}{n}}
={\lim\limits_{n\to\infty}}{\dfrac{F_{\alpha}^{n}({\bm x})}{n}}
$$
exists, and it only depends on $F_\alpha$.
We denote the fractional part of $\rho(F_{\alpha})$ (it means the fractional parts of all components of the vector $\rho(F_{\alpha})$) by $\rho(f_{\alpha})$ which is called the rotation vector of $f_\alpha$.
\end{definition}

\begin{remark}
In fact, for the general self-homeomorphisms (not rotations) of $T^n$, we can not define the rotation vectors of them by the same way,
because the limit
$$
{\lim\limits_{n\to\infty}}{\dfrac{F^n({\bm x})-{\bm x}}{n}}
={\lim\limits_{n\to\infty}}{\dfrac{F^n({\bm x})}{n}}
$$
may be dependent on the choices of ${\bm x}$, or does not exist.
\end{remark}

Assume that $f,\,g \in \mathscr{M}_{T^n}$ are two rotations of $T^n$.
Then $f$ and $g$ are topologically conjugate if and only if
$$
\rho(g)={\bm A}\rho(f)\quad({\rm mod}\,\,\,\mathbb{Z}^n),
$$
where the matrix ${\bm A}\in{\rm GL}_n(\mathbb {Z})$ is just the matrix form of the isomorphism $h_*:\,\pi_1(T^n) \rightarrow \pi_1(T^n)$ induced by $h$ which is the topological conjugacy from $f$ to $g$.
This fact implies that we can give a complete topologically conjugate classification of the rotations of $T^n$ by rotation vectors.

Next, we review some definitions and results about lens spaces \cite{HA-2000} and Lie groups \cite{PJF-1977,SMR-2007}.

\begin{definition}\label{def:01}
Regard $S^3$ as the unit sphere in $\mathbb{C}^2$, i.e.,
$$
S^3=\{(z_1,\,z_2);\,\,|z_1|^2+|z_2|^2=1,\,\,z_1,\,z_2 \in \mathbb{C}\}.
$$
Define a period homeomorphism $f:\,S^3 \rightarrow S^3$ by
$$
f:\,(z_1,\,z_2) \mapsto ({\textrm{e}}^{2\pi{\rm i}/p}z_1,\,{\textrm{e}}^{2\pi{\rm i} q/p}z_2),\qquad
\forall\,(z_1,\,z_2) \in S^3,
$$
where $p,\,q\in \mathbb{Z}$ and ${\rm gcd}\,(p, q)=1$.
Then we denote the quotient space $S^3/f$ by $L(p, q)$ which is called lens space.
\end{definition}

For any compact connected Lie group $G$, there exists some maximal commutative subgroup homeomorphic to $T^n$ which is called the maximal torus of $G$. We denote the maximal torus of $G$ by $T_G$.
The maximal torus of any compact connected Lie group may be not unique.
And a compact connected Lie group $G$ is commutative if and only if
$$
G=T_G \cong T^n.
$$

Fix one maximal torus $T_G$ of $G$. Then for any $g \in G$, there exist some $t \in T_G, s \in G$ such that $sg=ts$.
Assume that $g,\,s \in G$ and $t \in T_G$ satisfy $sg=ts$. Then we study the left actions $L_g,\,L_t$, and $L_s$.
For any $g' \in G$, we have $sgg'=tsg'$, and then
$$
L_s \circ L_g(g')=L_t \circ L_s(g'),
\quad\hbox{i.e.}\quad
L_s \circ L_g=L_t \circ L_s.
$$
Obviously, every left action is a self-homeomorphism of $G$.
Hence we have the following results.

Fix one maximal torus $T_G$ of $G$. Then any left action $L_g$ induced by an element $g$ in $G$ can be topologically conjugate to some left action $L_t$ induced by an element $t$ in $T_G$.

Therefore, it suffices to consider the left actions induced by the elements in the set $\mathscr{M}_{T_G}=\{L_g:\,G\rightarrow G;\,\,g \in T_G\}$ to classify all left actions on $G$.
And we have the following definition.

\begin{definition}\label{def:02}
Assume that $T_G$ is a maximal torus of a compact connected Lie group $G$, and $\Phi:\,T_G \rightarrow T^n$ is an isomorphism. Then $(T_G, \Phi)$ is called a maximal torus representation of $G$.
For every left action $L_g\in\mathscr{M}_{T_G}$, set $f=\Phi\circ L_g|_{T_G}\circ\Phi^{-1}$.
Since $T_G$ and $T^n$ are both commutative, we see that $f:\,T^n \rightarrow T^n$ is a rotation (left action) of the normal $n$-dimensional torus $T^n$. Then we define
$$
\rho(L_g)\triangleq\rho(f),
$$
and call $\rho(L_g)$ the rotation vector (number) of $L_g$ under the representation $(T_G, \Phi)$.
\end{definition}

Thus, it follows from Definition \ref{def:02} that we can define the rotation vectors of the left actions induced by the elements of the maximal tori of ${\rm SU}(2) \times T^n$,
and then utilize rotation vectors to give the complete topologically conjugate classifications of left actions.
Furthermore, we also discuss the relationship among their topological conjugacy, algebraic conjugacy and smooth conjugacy.

\section{Preliminaries}
In the first part of this section, we prove some useful lemmas.
In the second part, we define the rank and reductive rank of any set of real numbers $\{\alpha_i\}_{i=1}^n$,
and give three important applications of these two definitions under the background of the rotation vectors of the rotations of $T^n$.

\subsection{Some lemmas}
\begin{lemma}\label{lem:6}
Let $p_1,\,p_2,\,\cdots,\,p_n \in \mathbb{Z}$, and ${\rm gcd}\,(p_1, p_2, \cdots, p_n)=1$.
Then there exists some matrix ${\bm G} \in {\rm GL}_n(\mathbb{Z})$,
such that $(p_1, p_2,\,\cdots, p_n)$ is just some row of ${\bm G}$.
\end{lemma}

\begin{proof}
We can use the mathematical induction to prove this lemma.

When $n=1,\,2$, Lemma \ref{lem:6} is obviously true.
Suppose that if $n=m \geq 3$, Lemma \ref{lem:6} is true.
When $n=m+1$,
assume that ${\rm gcd}\,(p_2, p_3, \cdots, p_{m+1})=d$.
Then we have
$$
{\rm gcd}\,(\dfrac{p_2}{d}, \dfrac{p_3}{d}, \cdots, \dfrac{p_{m+1}}{d})=1,
\qquad {\rm gcd}\,(p_1, d)=1.
$$
It follows from the above assumption that there exists some matrix ${\bm G}' \in {\rm GL}_m(\mathbb{Z})$,
such that $(\dfrac{p_2}{d}, \dfrac{p_3}{d}, \cdots, \dfrac{p_{m+1}}{d})$ is just the $m$ row of ${\bm G}'$.
And the condition ${\rm gcd}\,(p_1, d)=1$ indicates that there exist some $s,\,t \in \mathbb{Z}$ such that $sp_1+td=1$.
Thus, take
$$
{\bm G}=\left(\begin{array}{cccc}
p_1&p_2&\cdots&p_{m+1}\\
{\bm 0}&&{\bm G}''&\\
-t&sp_2/d&\cdots&sp_{m+1}/d
\end{array}\right),
$$
where ${\bm G}''$ is the $(m-1) \times m$ matrix which is obtained by removing the $m$ row of ${\bm G}'$.
It is easy to verify that ${\bm G} \in {\rm GL}_{m+1}(\mathbb{Z})$.
Naturally, there exist some permutation matrices $\{{\bm T}_i\}_{i=1}^{m+1}$ such that $(p_1, p_2, \cdots, p_{m+1})$ is just the $i$ row of ${\bm T}_i{\bm G}$, where $i=1,\,2,\,\cdots,\,m+1$.
Obviously, ${\bm T}_i{\bm G} \in {\rm GL}_{m+1}(\mathbb{Z}),\,i=1, 2,\cdots, m+1$.
\end{proof}

\begin{lemma}\label{lem:1}\cite{HP-2018}
Suppose $G$ is a compact connected Lie group,
$L_g,\,L_{g'}$ are two left actions on $G$.
If $L_g$ and $L_{g'}$ are topologically conjugate,
then there exists some topological conjugacy $h$ from $L_g$ to $L_{g'}$ such that
$$
h(e)=e,
$$
where $e$ is the identity element of $G$.
\end{lemma}

\begin{lemma}\label{lem:2}
Suppose that $m \leq n$, $f_1$ is a rotation of $T^n$,
$f_2$ is a rotation of $T^m$,
$f:\,T^n \rightarrow T^m$ is a continuous subjection,
and
$$
f \circ f_1=f_2 \circ f,
$$
i.e., $f$ is a topological semi-conjugacy form $f_1$ to $f_2$.
For any
$$
{\bm k}=(\begin{array}{cccc}
k_1&k_2&\cdots&k_n
\end{array})^{\rm T} \in \pi_1(T^n),
$$
if $f_*:\,\pi_1(T^n) \rightarrow \pi_1(T^m)$ satisfies
$$
f_*:\,{\bm k} \mapsto {\bm A}{\bm k},\qquad {\bm A} \in {\rm{M}}_{m \times n}(\mathbb{Z}),
$$
then
$$
\rho(f_2)={\bm A}\rho(f_1)\quad(\rm{mod}\,\,\,\mathbb{Z}).
$$
\end{lemma}

\begin{proof}
Assume that $F_1:\,\mathbb{R}^n \rightarrow \mathbb{R}^n,\,F_2:\,\mathbb{R}^m \rightarrow \mathbb{R}^m$ are the lifts of $f_1$ and $f_2$ under the universal covering maps defined by
$$
F_1:\,{\bm x} \mapsto \bm{x}+\rho(f_1),\quad \forall\,{\bm x} \in \mathbb{R}^n,
\qquad F_2:\,{\bm x}' \mapsto {\bm x}'+\rho(f_2),\quad \forall\,{\bm x}' \in \mathbb{R}^m,
$$
respectively.
Since for any positive integer $n$, $\pi_1(\mathbb{R}^n)$ is trivial,
and $f \circ \pi$ is a continuous map,
then according to the map lifting theorem, there exists some continuous map $F:\,\mathbb{R}^n \rightarrow \mathbb{R}^m$, such that
$$
\pi' \circ F=f \circ \pi,
$$
where $\pi:\,\mathbb{R}^n \rightarrow T^n,\,\,\, \pi':\,\mathbb{R}^m \rightarrow T^m$ are universal covering maps.
Through a simple calculation, we have
$$
\pi' \circ (F \circ F_1)=(f \circ f_1) \circ \pi,\qquad \pi' \circ (F_2 \circ F)=(f_2 \circ f) \circ \pi.
$$
And together with the condition $f \circ f_1=f_2 \circ f$, one can see that
$$
\pi' \circ (F \circ F_1)=\pi' \circ (F_2 \circ F).
$$
Thus, by the definition of the universal covering map $\pi'$,
it is easy to verify that there exists some integer vector ${\bm l} \in \mathbb{Z}^m$, such that
$$
F_2 \circ F=F \circ F_1+{\bm l}.\eqno{(2.1)}
$$
Let $F_2$ act on the both sides of (2.1). We have
$$
F_2^2 \circ F=F_2 \circ (F \circ F_1+{\bm l}).
$$
And the definition of $F_2$ implies that
$$
F_2({\bm x}+{\bm l}')-F_2({\bm x})={\bm l}',
\qquad \forall\,{\bm x} \in \mathbb{R}^m,\,\,\,{\bm l}' \in \mathbb{Z}^m.
$$
So
\begin{align*}
F_2^2 \circ F
&=
F_2 \circ (F \circ F_1+{\bm l})
=F_2 \circ (F \circ F_1)+{\bm l}
\\
&=
(F_2 \circ F) \circ F_1+{\bm l}
=(F \circ F_1) \circ F_1+{\bm l}+{\bm l}
=F \circ F_1^2+2{\bm l}.
\end{align*}
Then let $F_2$ act on (2.1) by $n$ times, and we get
$$
F_2^n \circ F=F \circ F_1^n+n{\bm l},\qquad {\bm l} \in \mathbb{Z}^m,
$$
i.e.,
$$
\dfrac{F_2^n \circ F}{n}=\dfrac{F \circ F_1^n}{n}+{\bm l}.\eqno{(2.2)}
$$
Since for any
$$
{\bm k}=(\begin{array}{cccc}
k_1&k_2&\cdots&k_n
\end{array})^{\rm T} \in \pi_1(T^n)\,\cong\,\mathbb{Z}^n,
$$
the homomorphism $f_*:\,\pi_1(T^n) \rightarrow \pi_1(T^m)$ satisfies
$$
f_*:\,{\bm k} \mapsto {\bm A}{\bm k},\qquad {\bm A} \in {\rm{M}}_{m \times n}(\mathbb{Z}),
$$
then we can easily verify that
$$
F({\bm x}+{\bm k})-F({\bm x})={\bm A}{\bm k},\qquad \forall\,{\bm x} \in \mathbb{R}^n,\,\,\,{\bm k} \in \mathbb{Z}^n
$$
by the definition of the fundamental group $\pi_1(T^n)$.
Consequently, we obtain
$$
F(\bm{x}+\bm{k})-\bm{A}(\bm{x}+\bm{k})=F(\bm{x})-\bm{A}\bm{x}.\eqno{(2.3)}
$$
In fact, (2.3) indicates that $F(\bm{x})-\bm{A}\bm{x}$ is a continuous periodic function.
So there exists some $M > 0$, such that
$$
\|F(\bm{x})-\bm{A}\bm{x}\| \leq M,
$$
and then
$$
\|F(F_1^n(\bm{x}))-\bm{A}F_1^n(\bm{x})\| \leq M.
$$
Hence one can see that
$$
{\lim\limits_{n\to\infty}}\dfrac{\|F(F_1^n(\bm{x}))-\bm{A}F_1^n(\bm{x})\|}{n}=0,
$$
that means
$$
\lim\limits_{n\to\infty}\dfrac{F(F_1^n(\bm{x}))}{n}
=\lim\limits_{n\to\infty}\dfrac{\bm{A}F_1^n(\bm{x})}{n}.
$$
As $n\to\infty$, (2.2) implies that
$$
\lim\limits_{n\to\infty}\dfrac{F_2^n(F(\bm{x}))}{n}
=\lim\limits_{n\to\infty}\dfrac{F(F_1^n(\bm{x}))}{n}+{\bm l}
=\lim\limits_{n\to\infty}\dfrac{\bm{A}F_1^n(\bm{x})}{n}+{\bm l}.
$$
Therefore, it follows from Definition \ref{def:0} that
$$
\rho(F_2)={\bm A}\rho(F_1)+{\bm l}, \qquad {\bm l} \in \mathbb{Z}^m,
$$
and then
$$
{\rho (f_2)}={\bm A} {\rho (f_1)}\quad(\rm{mod}\,\,\,\mathbb{Z}).
$$

\end{proof}

\begin{lemma}\label{lem:5}\cite{HP-2018}
Assume that $G_1,\,G_2$ are topological groups,
$G_1',\,G_2'$ are the subgroups of $G_1$ and $G_2$, respectively,
and $h:\,G_1 \rightarrow G_2$ is a topological conjugacy from $\Gamma_{G_1'}$ to $\Gamma_{G_2'}$ satisfying $h(e_1)=e_2$,
where $\Gamma_{G_i'}$ is the $G_i'$ group action on $G_i$, and ${e_i}$ is the identity element of $G_i$, for $i=1,\,2$. Let
$\pi_1:\,G_1 \rightarrow G_1/G_1',\,\,\,\pi_2:\,G_2 \rightarrow G_2/G_2'$ be the quotient maps.
Then there exists some homeomorphism $h_\pi:\,G_1/G_1' \rightarrow G_2/G_2'$ induced by $h$ such that
$$
\pi_2 \circ h=h_\pi \circ \pi_1.
$$
In particular, if $G_1',\,G_2'$ are normal subgroups and $h$ is an isomorphism,
then $h_\pi$ is also an isomorphism induced by $h$.
\end{lemma}

\begin{remark}
In the present paper, for any topological group $G$ and its subgroup $G'$,
we always use the right coset of $G'$ to define the equivalence relation in $G$.
That means for any $g_1,\,g_2 \in G$, $[g_1]=[g_2]$ if and only if there exists some $g' \in G$ such that $g'g_1=g_2$, where $[g_1],\,[g_2]\in G/G'$ are the equivalence class of $g_1$ and $g_2$, respectively.
\end{remark}

\begin{lemma}\label{lem:3}
Define the composite map $f:\,M \rightarrow M$ by
$$
f=\pi \circ h \circ i:\,\xymatrixcolsep{1.5pc}
\xymatrix{
{\,\,M\,\,}  \ar[r]^-{i}
& {\,\,M \times T^n\,\,}  \ar[r]^-{h}
& {\,\,M \times T^n\,\,}  \ar[r]^-{\pi}
& {\,\,M\,\,},}
$$
where $n \geq 1$, $M$ is a $3$-dimensional lens space,
$i$ is the natural inclusion, $h$ is a self-homeomorphism on $M \times T^n$, $\pi$ is the projection.
Then $f$ is homotopic to some self-homeomorphism on $M$.
\end{lemma}

\begin{proof}
Since $h:\,M \times T^n \rightarrow M \times T^n$ is a homeomorphism,
and $\pi_1(T^n)\,\cong\,\mathbb{Z}^n$,
then one can see that $h_*|_{\pi_1(T^n)}$ is an automorphism of $\pi_1(T^n)$ satisfying $h_*|_{\pi_1(T^n)}={\bm A}$,
where $h_*$ is the automorphism of the fundamental group $\pi_1(M \times T^n)$ induced by $h$, ${\bm A} \in {\rm GL}_n(\mathbb{Z})$ is an $n \times n$ invertible integral matrix.
Thus, we can take some self-homeomorphism $H'$ on $M \times T^n$ satisfying
$$
H'_*|_{\pi_1(M)}={\rm id}_{\pi_1(M)},\qquad H'_*|_{\pi_1(T^n)}={\bm A}^{-1},
$$
where $H'_*$ is the automorphisms of the fundamental group $\pi_1(M \times T^n)$ induced by $H'$,
and ${\rm id}_{\pi_1(M)}$ is the identity endomorphism of $\pi_1(M)$.
Set $H=H' \circ h$.
Then we have
$$
H_*|_{\pi_1(T^n)}=(H'_*\circ h_*)|_{\pi_1(T^n)}=H'_*|_{\pi_1(T^n)} \circ h_*|_{\pi_1(T^n)}={\bm A}^{-1}{\bm A}={\rm id}_{\pi_1(T^n)},
$$
where ${\rm id}_{\pi_1(T^n)}$ is the identity endomorphism of $\pi_1(T^n)$.
Thus, by the definitions of $f$ and $H$ and the discussion in \cite{KR-2004} (see the section ``Proof of 3.1'' at page 8 to page 11),
it is obvious that
$$
f \times {\rm id}_{T^n}\,\simeq\,H.
$$

We regard $M \times T^{n-1} \times \mathbb{R}$ as the covering space of $M \times T^n\,\cong\,M \times T^{n-1} \times S^1$,
and define the covering map $\pi_0:\,M \times T^{n-1} \times \mathbb{R} \rightarrow M \times T^n\,\cong\,M \times T^{n-1} \times S^1$ by
$$
\pi_0 \triangleq {\rm id}_1 \times \pi_{\mathbb{R}},
$$
where ${\rm id}_1$ is the identity map of $M \times T^{n-1}$,
and $\pi_{\mathbb{R}}$ is the universal covering map from $\mathbb{R}$ to $S^1$.
Then investigate the following diagram.
$$
\xymatrixcolsep{2pc}
\xymatrix{
{\,\,M \times T^{n-1}\,\,} \ar[d]_-{{\rm id}_1} \ar[r]^-{i_1}
& {\,\,M \times T^{n-1} \times \mathbb{R}\,\,} \ar[d]_-{\pi_0} \ar[r]^-{h'}
& {\,\,M \times T^{n-1} \times \mathbb{R}\,\,} \ar[d]_-{\pi_0} \ar[r]^-{\pi_1}
& {\,\,M \times T^{n-1}\,\,} \ar[d]_-{{\rm id}_1}\\
{\,\,M \times T^{n-1}\,\,} \ar[r]_-{i'}
& {\,\,M \times T^{n-1} \times S^1 \,\,} \ar[r]_-{H}
& {\,\,M \times T^{n-1} \times S^1 \,\,} \ar[r]_-{\pi'}
& {\,\,M \times T^{n-1}\,\,}.}
$$
In the diagram, $H$ is the self-homeomorphism of $M \times T^n$ constructed above,
$h'$ is the lift of $H$ (according to the properties of $H$, the definition of $\pi_0$ and the map lifting theorem, we know that $H$ can be lifted under the covering map $\pi_0$),
$i_1,\,i'$ are the natural inclusion maps, $\pi_1,\,\pi'$ are the projections,
and $i_1,\,i',\,\pi_1,\,\pi',\,{\rm id}_1$ satisfy
$$
i_1(M \times T^{n-1})=M \times T^{n-1} \times \{0\},
\quad \pi_0 \circ i_1=i' \circ {\rm id}_1,
\quad {\rm id}_1 \circ \pi_1=\pi' \circ \pi_0.
$$
So this diagram is commutative,
and then we have
$$
f_1=\pi' \circ H \circ i'=\pi_1 \circ h' \circ i_1,
$$
In fact, $M \times T^{n-1}$ is a compact manifold, then obviously, there exists some $a \in \mathbb{R}$,
such that
$$
h'(M \times T^{n-1} \times [0, +\infty)) \subseteq M \times T^{n-1} \times (a, +\infty).
$$
Assume that
$$
W=M \times T^{n-1} \times [a, +\infty) \backslash h'(M \times T^{n-1} \times (0, +\infty)),
$$
and
$$
M_0=h'(M \times T^{n-1} \times \{0\}),\qquad M_1=M \times T^{n-1} \times \{a\}\,\cong\,M \times T^{n-1}.
$$
Then one can see that $(W;\,M_0,\,M_1)$ is a $h$-cobordism with $M_0$ and $M_1$.
Next, we prove $(W;\,M_0,\,M_1)$ is an $s$-cobordism satisfying $s$-cobordism theorem for three different cases.

(1) When $n \geq 3$,
using Whitehead torsion theory, and similar to the discussion in S. Kwasik and R. Schultz \cite{KR-2004}, we have
$$
\tau(f \times {\rm id}_{T^{n-1}})=\tau(W;\,M_0,\,M_1),
$$
where $\tau(f \times {\rm id}_{T^{n-1}})$ and $\tau(W;\,M_0,\,M_1)$ are the Whitehead torsions of $f \times {\rm id}_{T^{n-1}}$ and $(W;\,M_0,\,M_1)$, respectively.
And then according to the product formula of Whitehead torsion in M. M. Cohen \cite{MMC-1973}, we obtain
$$
\tau(f \times {\rm id}_{T^{n-1}})=\chi(T^{n-1})i_*\tau(f)+\chi(M)j_*\tau({\rm id}_{T^{n-1}})=0=\tau(W;\,M_0,\,M_1).
$$
So $h$-cobordism $(W;\,M_0,\,M_1)$ is an $s$-cobordism.
Since $n \geq 3$, then
$$
\dim(M_0)=\dim(M_1) \geq 5.
$$
Consequently, the investigations in \cite{HS-1969, KS-1977, JWM-1966} illuminate that $(W;\,M_0,\,M_1)$ satisfies $s$-cobordism theorem, that means
$$
M_1 \times [0, 1]\,\cong\,M \times T^{n-1} \times [0, 1]\,\cong\,W.
$$

(2) When $n=2$,
Similar to the discussion in (1), we know that $h$-cobordism $(W;\,M_0,\,M_1)$ is a $s$-cobordism.
In this case, one can see that $\dim(W)=5$, and
$$
\pi_1(W)\,\cong\,\pi_1(M_1)\,\cong\,\pi_1(M \times S^1)\,\cong\,\mathbb{Z}_p \oplus \mathbb{Z},
$$
where $p \in \mathbb{Z}_+$.
Then M. H. Freedman \cite{MHF-1983} indicates that $(W;\,M_0,\,M_1)$ satisfies $s$-cobordism theorem,
i.e.,
$$
M_1 \times [0, 1]\,\cong\,M \times S^1 \times [0, 1]\,\cong\,W.
$$

(3) When $n=1$, according to S. Kwasik and R. Schultz \cite{KR-1989}, we have
$$
M \times [0, 1]\,\cong\,W.
$$

In a word, for any positive integer $n$,  $(W;\,M_0,\,M_1)$ always satisfies
$$
M_1 \times [0, 1]\,\cong\,M \times T^{n-1} \times [0, 1]\,\cong\,W.
$$

Assume that $\mathcal{H}$ is a homeomorphism from $M \times T^{n-1} \times [0, 1]$ to $W$ satisfying
$$
\mathcal{H}(M \times T^{n-1} \times \{0\})=M_0,\qquad \mathcal{H}(M \times T^{n-1} \times \{1\})=M_1,
$$
and
$$
M_0'=M \times T^{n-1} \times \{0\}\,\cong\,M \times T^{n-1},
\qquad M_1'=M \times T^{n-1} \times \{1\}\,\cong\,M \times T^{n-1}.
$$
Take embedding maps
$$
f_0:\,M \times T^{n-1} \rightarrow M \times T^{n-1} \times [0, 1],
\qquad f'_0:\,M \times T^{n-1} \rightarrow M \times T^{n-1} \times [0, 1]
$$
satisfying
$$
f_0(M \times T^{n-1})=M'_0,\qquad f'_0(M \times T^{n-1})=M'_1.
$$
Then easy to see that $\mathcal{H} \circ f_0:\,M \times T^{n-1} \rightarrow W$ and $\mathcal{H} \circ f'_0:\,M \times T^{n-1} \rightarrow W$ are also embedding maps satisfying
$$
\mathcal{H} \circ f_0(M \times T^{n-1})=M_0,\qquad \mathcal{H} \circ f'_0(M \times T^{n-1})=M_1
$$
Define $F:\, M \times T^{n-1} \times [0, 1] \rightarrow M \times T^{n-1} \times [0, 1]$ by
$$
F:\,(x, t) \mapsto (h_0^{-1} \circ h'|_{M \times T^{n-1} \times \{0\}} \circ f_0(x),\,t),
\qquad \forall\,x \in M \times T^{n-1},\,t \in [0, 1],
$$
where $h_0=\mathcal{H}|_{M_0'} \circ f_0:\,M \times T^{n-1} \rightarrow M_0$ is a homeomorphism.
Next, we study the map $\mathcal{H} \circ F:\,M \times T^{n-1} \times [0, 1] \rightarrow W$.
One can see that for any $t \in [0, 1]$, $\mathcal{H} \circ F|_{M \times T^{n-1} \times \{t\}}$ is always an embedding map, and $\mathcal{H} \circ F$ satisfies
\begin{gather*}
\mathcal{H} \circ F(x, 0)=h'|_{M \times T^{n-1} \times \{0\}} \circ f_0(x),
\\
\mathcal{H} \circ F(x, 1)=h'_0 \circ h_0^{-1} \circ h'|_{M \times T^{n-1} \times \{0\}} \circ f_0(x),
\end{gather*}
where $h'_0=\mathcal{H}|_{M_1'} \circ f'_0:\,M \times T^{n-1} \rightarrow M_1$ is a homeomorphism.
Thus, $\mathcal{H} \circ F$ is an isotopy from $h'|_{M \times T^{n-1} \times \{0\}} \circ f_0$ to $h'_0 \circ h_0^{-1} \circ h'|_{M \times T^{n-1} \times \{0\}} \circ f_0$.
So
$$
\pi_1 \circ h'|_{M \times T^{n-1} \times \{0\}} \circ f_0 \simeq \pi_1 \circ h'_0 \circ h_0^{-1} \circ h'|_{M \times T^{n-1} \times \{0\}} \circ f_0,
$$
where $\pi_1:\,M \times T^{n-1} \times \mathbb{R} \rightarrow M \times T^{n-1}$ is the projection in the above diagram.
It is no hard to see that
$$
\pi_1 \circ h'|_{M \times T^{n-1} \times \{0\}} \circ f_0=\pi_1 \circ h' \circ i_1=f_1,
$$
and
$$\pi_1 \circ h'_0 \circ h_0^{-1} \circ h'|_{M \times T^{n-1} \times \{0\}} \circ f_0
$$
is a self-homeomorphism on $M \times T^{n-1}$.
Set
$$
h_1=\pi_1 \circ h'_0 \circ h_0^{-1} \circ h'|_{M \times T^{n-1} \times \{0\}} \circ f_0.
$$
We know that
$$
f_1=\pi_1 \circ h' \circ i_1=\pi' \circ H \circ i'.
$$
Then the construction of $H$ indicates $f_1 \times {\rm id}_{S^1}\,\simeq\,H$.
And combining this condition with the fact $f \times {\rm id}_{T^n}\,\simeq\,H$,
we obtain
$$
f \times {\rm id}_{T^{n-1}}\,\simeq\,f_1\,\simeq\,h_1.
$$

Next, we regard $M \times T^{n-2} \times \mathbb{R}$ as the covering space $M \times T^{n-1}\,\cong\,M \times T^{n-2} \times S^1$,
and define the covering map $\pi_0':\,M \times T^{n-2} \times \mathbb{R} \rightarrow M \times T^n\,\cong\,M \times T^{n-2} \times S^1$ by
$$
\pi_0' \triangleq {\rm id}_2 \times \pi_{\mathbb{R}},
$$
where ${\rm id}_2$ is the identity map of $M \times T^{n-2}$,
and $\pi_{\mathbb{R}}$ is the universal covering map from $\mathbb{R}$ to $S^1$.
Then investigate the following diagram.
$$
\xymatrixcolsep{2pc}
\xymatrix{
{\,\,M \times T^{n-2}\,\,} \ar[d]_-{{\rm id}_2} \ar[r]^-{i_2}
& {\,\,M \times T^{n-2} \times \mathbb{R}\,\,} \ar[d]_-{\pi'_0} \ar[r]^-{h''}
& {\,\,M \times T^{n-2} \times \mathbb{R}\,\,} \ar[d]_-{\pi'_0} \ar[r]^-{\pi_2}
& {\,\,M \times T^{n-2}\,\,} \ar[d]_-{{\rm id}_2}\\
{\,\,M \times T^{n-2}\,\,} \ar[r]_-{i''}
& {\,\,M \times T^{n-2} \times S^1 \,\,} \ar[r]_-{h_1}
& {\,\,M \times T^{n-2} \times S^1 \,\,} \ar[r]_-{\pi''}
& {\,\,M \times T^{n-2}\,\,}.}
$$
Set
$$
f_2=\pi_2 \circ h'' \circ i_2=\pi'' \circ h_1 \circ i''.
$$
By the same way, one can see that there exists some self-homeomorphism $h_2$ on $M \times T^{n-2}$ such that $f_2\,\simeq\,h_2$,
and then we have
$$
f \times {\rm id}_{T^{n-2}}\,\simeq\,f_2\,\simeq\,h_2.
$$
So we proceed to do the same work by $n-2$ times.
Finally, we can find a self-homeomorphism $h_n$ on $M$ such that
$$
f\,\simeq\,f_n\,\simeq\,h_n.
$$

\end{proof}

\begin{lemma}\label{lem:4}\cite{SH-2010}
Suppose that $f$ is a degree 1 self-map on $L(p, q)$.
Then $f$ is homotopic to an orientation-preserving homeomorphism if and only if
$$
f_*(l)=\left\{
\begin{array}{l}
\pm\,l,\qquad\qquad p \nmid q^2-1,
\\[0.25 cm]
\pm\,l, \pm ql,\qquad p \mid q^2-1,
\end{array}
\right.
\qquad \forall\,l \in \pi_1(L(p, q)),
$$
where $f_*$ is the endomorphism of the fundamental group $\pi_1(L(p, q))$ induced by $f$.
\end{lemma}

\begin{lemma}\label{lem:7}\cite{HP-2018}
Assume that $B$ is a topological space,
$(E, \pi)$ is the covering space of $B$,
$f,\,g$ are two continuous self-maps of $B$,
$h$ is a self-homeomorphism of $B$, and $f,\,g,\,h$ can be lifted under the covering map $\pi$.
Then $f$ and $g$ are topologically conjugate,
and the homeomorphism $h:\,B \rightarrow B$ is a topological conjugacy from $f$ to $g$, that means
$$
h \circ f=g \circ h
$$
if and only if for any fixed lift $\tilde{f}$ of $f$ and any fixed lift $\tilde{h}$ of $h$,
there exists some certain lift $\tilde{g}$ of $g$ such that $\tilde{f}$ and $\tilde{g}$ are topologically conjugate,
and the homeomorphism $\tilde{h}:\,E \rightarrow E$ is a topological conjugacy from $\tilde{f}$ to $\tilde{g}$, that means
$$
\tilde{h} \circ \tilde{f}= \tilde{g} \circ \tilde{h}.
$$
\end{lemma}

\begin{lemma}\label{lem:8}
Let $k_0,\,k_1,\cdots,k_m \in \mathbb{Z}$ with ${\rm gcd}\,(k_0, k_1, \cdots, k_m)=1$,
and
$$
G=\left\{\left(\begin{array}{cc}
z_0z_1 \cdots z_n&0\\
0&\bar{z_0}
\end{array}\right);\,z_0^{k_0}z_1^{k_1}z_2^{k_2} \cdots z_m^{k_m}=1,\, z_i=1,\,i=m+1,\,\cdots,\,n \right\},
$$
where $z_j \in \mathbb{C},\,|z_j|=1,\,j=0,\,1,\,\cdots,\,m$.
Then $G\,\cong\,T^m$ is a subgroup of ${\rm SU}(2) \times T^n$, and
$$
{\rm SU}(2) \times T^n/G\,\cong\,L(k_0, -1) \times T^{n-m},
$$
where $L(p, -1)$ is a $3$-dimensional lens space.
\end{lemma}

\begin{proof}
In fact, It is easy to verify that $G$ is a subgroup of ${\rm SU}(2) \times T^n$.
And together with the condition
$$
{\rm gcd}\,(k_0, k_1, \cdots, k_m)=1,
$$
we have $G\,\cong\,T^m$.

Regard $S^3$ as the unit sphere in quaternion space, and then denote the elements of $S^3$ by $z+{\rm j}z'$, where $|z|^2+|z'|^2=1,\,\,z,\,z' \in \mathbb{C}$.
Notice that ${\rm SU}(2) \cong S^3$, and $S^3$ can be obtained by gluing the outside surfaces of two solid tori denoted by $D^2_+ \times S^1$ and $D^2_- \times S^1$, where $D^2_+,\,D^2_-$ are two unit closed disks.
Define continuous maps $f_1:\,D^2_+ \times S^1 \rightarrow S^3,
\,f_2:\,D^2_- \times S^1 \rightarrow S^3$ by
\begin{gather*}
f_1:\,(\lambda, z_0)\mapsto\dfrac{z_0(1+{\rm j}\lambda)}{\sqrt{|\lambda|^2+1}},
\qquad \forall\, (\lambda, z_0) \in D^2_+ \times S^1,
\\
f_2:\,(\lambda', z_0')\mapsto\dfrac{z_0'(\lambda'+{\rm j})}{\sqrt{|\lambda'|^2+1}},
\qquad \forall\, (\lambda', z_0') \in D^2_- \times S^1,
\end{gather*}
respectively.
Set
$$
f=f_1 \sqcup f_2:\,D^2_+ \times S^1 \sqcup D^2_- \times S^1 \rightarrow S^3,
$$
where ``$\sqcup$'' means disjoint union.
Obviously, $f$ is a surjection.
So $f$ can induce an equivalence relation ``$\sim$'' defined by
$$
(\lambda, z_0) \sim (\lambda', z_0') \Longleftrightarrow
f_1(\lambda, z_0)=f_2(\lambda', z_0').
$$
Then
$$
D^2_+ \times S^1 \sqcup D^2_- \times S^1/\sim\,\,\,\cong S^3.
$$
According to the equation $f_1(\lambda, z_0)=f_2(\lambda', z'_0)$, we get
$$
\left\{
\begin{array}{l}
\dfrac{z_0}{\sqrt{|\lambda|^2+1}}=\dfrac{z'_0\lambda'}{\sqrt{|\lambda'|^2+1}},
\\[0.5 cm]
\dfrac{z_0\bar{\lambda}}{\sqrt{|\lambda|^2+1}}=\dfrac{z'_0}{\sqrt{|\lambda'|^2+1}},
\end{array}
\right.
\quad\hbox{i.e.,}\quad
\left\{
\begin{array}{l}
\lambda'=\lambda,
\\[0.25 cm]
z_0'=z_0\bar{\lambda}.
\end{array}
\right.
$$
Therefore, we can define an identification map $F:\,\partial D^2_+ \times S^1 \rightarrow \partial D^2_- \times S^1$ by
$$
F:\,(\lambda, z_0)\mapsto(\lambda', z'_0)=(\lambda, z_0\bar{\lambda}),
\qquad \forall\,(\lambda, z_0) \in \partial D^2_+ \times S^1.
$$
By this way, we regard $S^3$ as a circle bundle on $S^2$.
If we denote the arguments of $\lambda_1,\,\lambda_2,\,z_1,\,z_2$ by $\gamma,\,\gamma',\,\alpha_0,\,\alpha'_0$, respectively,
then $F$ corresponds to a $2\times2$ matrix which describes the relationship of arguments as follows
$$
\left(\arraycolsep=1pt\begin{array}{c}
\gamma'\\
\alpha_0'
\end{array}\right)=\left(\begin{array}{cc}
1&0\\
-1&1
\end{array}\right)\left(\arraycolsep=1pt\begin{array}{c}
\gamma\\
\alpha_0\
\end{array}\right).
$$
Naturally, ${\rm SU}(2) \times T^n$can be regarded as a $T^{n+1}$ fibre bundle on $S^2$,
and every fiber is a maximal torus of ${\rm SU}(2) \times T^n$.
Here, we state that the symbol ``$\times$'' between two Lie groups in the present paper just denotes Cartesian product in topology, not direct sum in algebra.
Then the identification map $F$ can induce an identification map $F':\,\partial D^2_+ \times S^1 \times T^n \rightarrow \partial D^2_- \times S^1 \times T^n$ in ${\rm SU}(2) \times T^n$ satisfying for any $(\lambda, z_0, z_1, \cdots, z_n) \in \partial D^2_+ \times S^1 \times T^n$,
$$
F':\,(\lambda, z_0, z_1, \cdots, z_n) \mapsto (\lambda, z_0\bar{\lambda}, z_1, \cdots, z_n).
$$
Denote the arguments of $z_i,\,z'_i$ by $\alpha_i,\,\alpha'_i$, where $i=1,\,2,\,\cdots,\,n$.
Then $F'$ corresponds to a square matrix with order $n+2$ which describes the relationship of arguments, i.e.,
$$
\left(\arraycolsep=1pt\begin{array}{c}
\gamma'\\
\alpha'_0\\
\vdots\\
\alpha'_n
\end{array}\right)=\left(\begin{array}{ccc}
1&0&{\bm 0}\\
-1&1&{\bm 0}\\
{\bm 0}&{\bm 0}&{\bm I}_n
\end{array}\right)\left(\arraycolsep=1pt\begin{array}{c}
\gamma\\
\alpha_0\\
\vdots\\
\alpha_n
\end{array}\right),
$$
where ${\bm I}_n$ is the identity matrix with order $n$.
Since ${\rm gcd}\,(k_0, k_1, \cdots, k_m)=1$,
and then Lemma \ref{lem:6} indicates that there exists some matrix ${\bm G} \in {\rm GL}_{m+1}(\mathbb{Z})$ such that $(k_0, k_1, \cdots, k_m)$ is just the first row of ${\bm G}$.
Through a transformation of arguments, we have
$$
\left(\arraycolsep=1pt\begin{array}{c}
\gamma\\
\beta_0\\
\vdots\\
\beta_m\\
\alpha_{m+1}\\
\vdots\\
\alpha_n
\end{array}\right)={\bm A}\left(\arraycolsep=1pt\begin{array}{c}
\gamma\\
\alpha_0\\
\vdots\\
\alpha_m\\
\alpha_{m+1}\\
\vdots\\
\alpha_n
\end{array}\right),\,\,\,
\left(\arraycolsep=1pt\begin{array}{c}
\gamma'\\
\beta'_0\\
\vdots\\
\beta'_m\\
\alpha'_{m+1}\\
\vdots\\
\alpha'_n
\end{array}\right)={\bm A}\left(\arraycolsep=1pt\begin{array}{c}
\gamma'\\
\alpha'_0\\
\vdots\\
\alpha'_m\\
\alpha'_{m+1}\\
\vdots\\
\alpha'_n
\end{array}\right),
$$
where ${\bm A}=\left(\begin{array}{ccc}
1&{\bm 0}&{\bm 0}\\
{\bm 0}&{\bm G}&{\bm 0}\\
{\bm 0}&{\bm 0}&{\bm I}_{n-m}
\end{array}\right)$,
and ${\bm I}_{n-m}$ is the identity matrix with order $n-m$.
Obviously, $|{\bm A}|=1$.
Then
$$
\left(\arraycolsep=1pt\begin{array}{c}
\gamma\\
\alpha_0\\
\vdots\\
\alpha_m\\
\alpha_{m+1}\\
\vdots\\
\alpha_n
\end{array}\right)={\bm A}^{-1}\left(\arraycolsep=1pt\begin{array}{c}
\gamma\\
\beta_0\\
\vdots\\
\beta_m\\
\alpha_{m+1}\\
\vdots\\
\alpha_n
\end{array}\right),\qquad
\left(\arraycolsep=1pt\begin{array}{c}
\gamma'\\
\alpha'_0\\
\vdots\\
\alpha'_m\\
\alpha'_{m+1}\\
\vdots\\
\alpha'_n
\end{array}\right)={\bm A}^{-1}\left(\arraycolsep=1pt\begin{array}{c}
\gamma'\\
\beta'_0\\
\vdots\\
\beta'_m\\
\alpha'_{m+1}\\
\vdots\\
\alpha'_n
\end{array}\right).
$$
Thus, we have
$$
\left(\arraycolsep=1pt\begin{array}{c}
\gamma'\\
\beta'_0\\
\vdots\\
\beta'_m\\
\alpha'_{m+1}\\
\vdots\\
\alpha'_n
\end{array}\right)={\bm A}\left(\begin{array}{ccc}
1&0&{\bm 0}\\
-1&1&{\bm 0}\\
{\bm 0}&{\bm 0}&{\bm I}_n
\end{array}\right){\bm A}^{-1}\left(\arraycolsep=1pt\begin{array}{c}
\gamma\\
\beta_0\\
\vdots\\
\beta_m\\
\alpha_{m+1}\\
\vdots\\
\alpha_n
\end{array}\right),
$$
It follows from a simple calculation that we get
$$
\left(\arraycolsep=1pt\begin{array}{c}
\gamma'\\
\beta'_0\\
\vdots\\
\beta'_m\\
\alpha'_{m+1}\\
\vdots\\
\alpha'_n
\end{array}\right)=\left(\begin{array}{cc}
{\bm G}'&{\bm 0}\\
{\bm 0}&{\bm I}_{n-m}
\end{array}\right)\left(\arraycolsep=1pt\begin{array}{c}
\gamma\\
\beta_0\\
\vdots\\
\beta_m\\
\alpha_{m+1}\\
\vdots\\
\alpha_n
\end{array}\right),\eqno{(2.4)}
$$
where ${\bm G}'$ is a lower triangular matrix with order $m+2$ whose diagonal elements are all $1$,
and the element lying in the second row and the fist column of ${\bm G}'$ is $-k_0$.
In fact, $m$-dimensional subgroup $G\,\cong\,T^m$ is included in some maximal torus of ${\rm SU}(2) \times T^n$.
Then one can see that the quotient space ${\rm SU}(2) \times T^n/G$ is a $T^{n-m+1}$ fibre bundle on $S^2$.
Thus, the identification map $F'$ naturally induces an identification map
$$
F'':\,\partial D^2_+ \times S^1 \times T^{n-m}\rightarrow \partial D^2_- \times S^1 \times T^{n-m}
$$
in ${\rm SU}(2) \times T^n/G$.
According to (2.4), we see that $F''$ corresponds to a square matrix of order $n-m+2$ which describes the relationship of arguments.
And this matrix is just the matrix in (2.4) except the $3,\,4,\,\cdots,\,m+2$ row and the $3,\,4,\,\cdots,\,m+2$ column, i.e.,
$$
\left(\arraycolsep=1pt\begin{array}{c}
\gamma'\\
\beta'_0\\
\alpha'_{m+1}\\
\vdots\\
\alpha'_n
\end{array}\right)=\left(\begin{array}{ccc}
1&0&{\bm 0}\\
-k_0&1&{\bm 0}\\
{\bm 0}&{\bm 0}&{\bm I}_{n-m}
\end{array}\right)\left(\arraycolsep=1pt\begin{array}{c}
\gamma'\\
\beta_0\\
\alpha_{m+1}\\
\vdots\\
\alpha_n
\end{array}\right).
$$
Then the identification map $F''$ satisfies
$$
F''(\lambda, z_0, z_{m+1}, \cdots z_{n})=(\lambda, z_0\bar{\lambda}^{k_0}, z_{m+1}, \cdots z_{n}),
$$
for any $(\lambda, z_0, z_{m+1}, \cdots z_{n})\in\partial D^2_+ \times S^1 \times T^{n-m}$.
Using the above observations, we have
\begin{gather*}
{\rm SU}(2) \times T^n \cong D^2_+ \times S^1 \times T^n \sqcup D^2_- \times S^1 \times T^n/\sim_1,
\\
{\rm SU}(2) \times T^n/G \cong D^2_+ \times S^1 \times T^{n-m}\sqcup D^2_- \times S^1 \times T^{n-m}/\sim_2,
\end{gather*}
where ``$\sim_1$'' and ``$\sim_2''$ are the equivalence relations induced by $F'$ and $F''$, respectively.
We can denote the elements of ${\rm SU}(2) \times T^n$ and ${\rm SU}(2) \times T^n/G$ by $(z, {\textrm{e}}^{2\pi{\rm i}\alpha_0}, \cdots, {\textrm{e}}^{2\pi{\rm i}\alpha_n})$ and $(z, {\textrm{e}}^{2\pi{\rm i}\theta_0}, {\textrm{e}}^{2\pi{\rm i}\theta_1},\cdots, {\textrm{e}}^{2\pi{\rm i}\theta_{n-m}})$, respectively,
where $|z| \leq 1$, and $z \in D^2_+$ or $ D^2_-$.
Thus, it is not difficult to prove that the quotient map $\pi:\,{\rm SU}(2) \times T^n \rightarrow {\rm SU}(2) \times T^n/G$ satisfies
$$
\pi(z, {\textrm{e}}^{2\pi{\rm i}\alpha_0}, \cdots, {\textrm{e}}^{2\pi{\rm i}\alpha_n})=(z, {\textrm{e}}^{2\pi{\rm i}\beta_0}, {\textrm{e}}^{2\pi{\rm i}\alpha_{m+1}},\cdots, {\textrm{e}}^{2\pi{\rm i}\alpha_n}),
$$
for any $(z, {\textrm{e}}^{2\pi{\rm i}\alpha_0}, \cdots, {\textrm{e}}^{2\pi{\rm i}\alpha_n})\in{\rm SU}(2) \times T^n$,
where $\beta_0=\sum_{i=0}^mk_i\alpha_i$.
On the other hand, define two maps $\pi_+:D^2_+ \times S^1 \times T^{n-m} \rightarrow D^2_+ \times S^1 \times T^{n-m},
\,\,\,\pi_-:D^2_- \times S^1 \times T^{n-m} \rightarrow D^2_- \times S^1 \times T^{n-m}$ by
\begin{gather*}
\pi_+:\,(\lambda, z_0, z_1, \cdots, z_{n-m}) \mapsto (\lambda, z_0^{k_0}, z_1, \cdots, z_{n-m}),
\\
\pi_-:\,(\lambda', z'_0, z'_1, \cdots, z'_{n-m}) \mapsto (\lambda', z_0'^{k_0}, z'_1, \cdots, z'_{n-m}),
\end{gather*}
respectively.
And through a simple verification, we have
$$
\pi_- \circ (F \times {\rm id}_{T^{n-m}})=F'' \circ \pi_+,
$$
where ${\rm id}_{T^{n-m}}$ is the identity map of $T^{n-m}$,
$$
F \times {\rm id}_{T^{n-m}}:\,D^2_+ \times S^1 \times T^{n-m} \rightarrow D^2_- \times S^1 \times T^{n-m}
$$
is an identification map, and the equivalence relation ``$\sim_3$" induced by $F \times {\rm id}_{T^{n-m}}$ satisfies
$$
{\rm SU}(2) \times T^{n-m} \cong D^2_+ \times S^1 \times T^{n-m} \sqcup D^2_- \times S^1 \times T^{n-m}/\sim_3.
$$
This fact illustrates that ${\rm SU}(2) \times T^{n-m}$ is the $k_0$-fold covering space of ${\rm SU}(2) \times T^n/G$, and the local representation of the covering map from ${\rm SU}(2) \times T^{n-m}$ to ${\rm SU}(2) \times T^n/G$ is $\{\pi_+,\,\pi_-\}$.
One can see that ${\rm SU}(2) \times T^{n-m}$ is the $k_0$-fold covering space of $L(k_0, -1) \times T^{n-m}$.
And according to one definition of $L(k_0, -1)$, the local representation of the covering map from ${\rm SU}(2) \times T^{n-m}$ to $L(k_0, -1) \times T^{n-m}$ is just $\{\pi_+,\,\pi_-\}$.
Therefore, we obtain
$$
{\rm SU}(2) \times T^n/G \cong L(k_0, -1) \times T^{n-m}.
$$

\end{proof}

\subsection{Definitions and applications of rank and reductive rank}

For any set of real numbers $\{\alpha_i\}_{i=1}^n$, we give the following definition.

\begin{definition}\label{rrr}
Regard $\mathbb{R}$ as an infinite dimensional vector space over the rational field $\mathbb{Q}$.
Let $\{\alpha_i\}_{i=1}^n$ be $n$ points in vector space $\mathbb{R}$.
Define the rank of $\{\alpha_i\}_{i=1}^n$ over $\mathbb{Q}$ by
\begin{align*}
\mathcal{R}(\alpha_1, \alpha_2, \cdots, \alpha_n)
&={\rm{dim}}({\rm{Span}}_{\mathbb{Q}}\{\alpha_i\}_{i=1}^n)
\\
&\triangleq {\rm{dim}}(\{\sum^n_{i=1}r_i\alpha_i;\,r_i \in \mathbb{Q},\,i=1, 2, \cdots, n\}).
\end{align*}
Moreover, define the reductive rank of $\{\alpha_i\}_{i=1}^n$ over $\mathbb{Q}$ by
$$
\widetilde{\mathcal{R}}(\alpha_1, \alpha_2, \cdots, \alpha_n)={\rm{dim}}({\rm{Span}}_{\mathbb{Q}}\{1, \alpha_1, \alpha_2, \cdots, \alpha_n\}/\mathbb{Q}).
$$
\end{definition}

\begin{remark}\label{rem:1}
Let $\{\alpha_i\}_{i=1}^n$ be $n$ points in vector space $\mathbb{R}$ over the rational field $\mathbb{Q}$.
Then there exists a set of linearly independent real numbers $\{\alpha_{k_j}\}_{j=1}^m \subseteq \{\alpha_i\}_{i=1}^n$,
i.e., if $\sum^{m}_{j=1}r_{k_j}\alpha_{k_j}=0$, then $r_{k_j}=0$,
where $r_{k_j} \in \mathbb{Q},\,j=1,\,2,\,\cdots,\,m$.
Assume that $\{\alpha_{k_j}\}_{j=1}^{m}$ is a maximal linearly independent subset of $\{\alpha_i\}_{i=1}^n$.
Then according to Definition \ref{rrr}, we know
$$
\mathcal{R}(\alpha_1, \alpha_2, \cdots, \alpha_n)=m.
$$
For reductive rank, similarly, there exists a set of rationally independent real numbers $\{\alpha_{l_{j'}}\}_{j'=1}^{m'} \subseteq \{\alpha_i\}_{i=1}^n$,
i.e., if $\sum^{m'}_{j'=1}r_{l_{j'}}\alpha_{l_{j'}}+r=0$, then $r=r_{l_{j'}}=0$,
where $r,\,r_{k_{j'}} \in \mathbb{Q},\,j'=1,\,2,\,\cdots,\,m$.
If $\{\alpha_{l_{j'}}\}_{j'=1}^{m'}$ is a maximal rationally independent subset of $\{\alpha_i\}_{i=1}^n$,
then Definition \ref{rrr} implies that
$$
\widetilde{\mathcal{R}}(\alpha_1, \alpha_2, \cdots, \alpha_n)=m'.
$$
Furthermore, together the above discussion with Definition \ref{rrr}, it is not difficult to see that
$$
\widetilde{\mathcal{R}}(\alpha_1, \alpha_2, \cdots, \alpha_n)=\mathcal{R}(\alpha_1, \alpha_2, \cdots, \alpha_n)
$$
or
$$
\widetilde{\mathcal{R}}(\alpha_1, \alpha_2, \cdots, \alpha_n)=\mathcal{R}(\alpha_1, \alpha_2, \cdots, \alpha_n)-1.
$$
If $\widetilde{\mathcal{R}}(\alpha_1, \alpha_2, \cdots, \alpha_n)=\mathcal{R}(\alpha_1, \alpha_2, \cdots, \alpha_n)=m$,
we know that the maximal linearly independent subset $\{\alpha_{k_j}\}_{j=1}^{m}$ of $\{\alpha_i\}_{i=1}^n$ is also a rationally independent subset,
and $\alpha_{k_1},\,\alpha_{k_2},\,\cdots,\,\alpha_{k_m}$ are $m$ rationally independent irrational numbers.
When $\widetilde{\mathcal{R}}(\alpha_1, \alpha_2, \cdots, \alpha_n)=\mathcal{R}(\alpha_1, \alpha_2, \cdots, \alpha_n)-1$.
Let $\{\alpha_{k_j}\}_{j=1}^{m}$ be a maximal linearly independent subset of $\{\alpha_i\}_{i=1}^n$.
And assume that $\{\alpha_{k_j}\}_{j=1}^{m-1}$ is a maximal rationally independent subset of $\{\alpha_i\}_{i=1}^n$.
Then, in this case, $\alpha_{k_1}, \alpha_{k_2}, \cdots, \alpha_{k_m}$ are rationally dependent,
$\alpha_{k_1}, \alpha_{k_2}, \cdots, \alpha_{k_{m-1}}$ are $m-1$ rationally independent irrational numbers,
and $\alpha_{k_m}$ is either a rational number or an irrational number.
\end{remark}

Since $\mathbb{R}$ can be regarded as an infinite dimensional vector space over the rational field $\mathbb{Q}$,
so $\{\alpha_i\}_{i=1}^n$ is naturally a subspace of $\mathbb{R}$.
Thus, according to the above analysis and the definition of the dimension of vector space in linear algebra theory,
it is easy to verify that the definitions of rank and reductive rank are well-defined.

Next, we use the above properties to give three important applications of rank and reductive rank under the background of the rotations of $T^n$.
First of all, for a rotation $f$ of $T^n$ with $\rho(f)=\left(\begin{array}{cccc}
\theta_1&\theta_2&\cdots&\theta_n
\end{array}\right)^{\rm T}$,
where $\theta_1,\,\theta_2,\,\cdots,\,\theta_n \in [0, 1)$,
we define
$$
\mathcal{R}_{\rho(f)}\triangleq\mathcal{R}(\theta_1, \theta_2, \cdots, \theta_n),
\qquad\widetilde{\mathcal{R}}_{\rho(f)}\triangleq\widetilde{\mathcal{R}}(\theta_1, \theta_2, \cdots, \theta_n).
$$

\begin{proposition}\label{prop:1}
Assume that $\{\alpha_i\}^n_{i=1},\,\{\alpha'_j\}^n_{j=1} \subseteq \mathbb{R}$ satisfy
$$
\left(\begin{array}{c}
\alpha'_1\\
\alpha'_2\\
\vdots\\
\alpha'_n
\end{array}\right)={\bm A}\left(\begin{array}{c}
\alpha_1\\
\alpha_2\\
\vdots\\
\alpha_n
\end{array}\right),\qquad {\bm A} \in {\rm{GL}}_n(\mathbb{Z}).
$$
Then
$$
\mathcal{R}(\alpha_1, \alpha_2, \cdots, \alpha_n)=\mathcal{R}(\alpha'_1, \alpha'_2, \cdots, \alpha'_n).
$$
Furthermore, if $f,\,g$ are two rotations on $T^n$ satisfying
$$
\rho(f)=\left(\begin{array}{c}
\alpha_1\\
\alpha_2\\
\vdots\\
\alpha_n
\end{array}\right)\,\,\,({\rm{mod}}\,\,\,\mathbb{Z}),\qquad \rho(g)=\left(\begin{array}{c}
\alpha'_1\\
\alpha'_2\\
\vdots\\
\alpha'_n
\end{array}\right)\,\,\,({\rm{mod}}\,\,\,\mathbb{Z}),
$$
and
$$
\rho(g)={\bm A}'\rho(f)\quad({\rm{mod}}\,\,\,\mathbb{Z}),
$$
where ${\bm A}' \in {\rm{GL}}_n(\mathbb{Z})$,
then
$$
\widetilde{\mathcal{R}}_{\rho(f)}=\widetilde{\mathcal{R}}_{\rho(g)}.
$$
\end{proposition}

\begin{proof}
By the definitions of rank and reductive rank, it is not difficult to prove this proposition.
\end{proof}

\begin{proposition}\label{prop:2}
Let $f$ be a rotation of $T^n$ satisfying
$$
\rho(f)=(\begin{array}{cccc}
\alpha_1&\alpha_2&\cdots&\alpha_n
\end{array})^{\rm T},\qquad \alpha_i \in [0, 1),\,\,\,i=1, 2, \cdots, n.
$$
If $\mathcal{R}_{\rho(f)}=m$,
then there exists some rotation $g$ of $T^n$ satisfying
$$
\rho(g)=(\begin{array}{ccccccc}
\beta_1&\beta_2&\cdots&\beta_m&0&\cdots&0
\end{array})^{\rm T},\qquad \beta_1,\,\beta_2,\,\cdots,\,\beta_m \in [0, 1),
$$
such that $f$ and $g$ are topologically conjugate.
Furthermore, if $\widetilde{\mathcal{R}}_{\rho(f)}=\mathcal{R}_{\rho(f)}=m$,
then $\beta_1,\,\beta_2,\,\cdots,\,\beta_m$ are $m$ rationally independent irrational numbers.
If $\widetilde{\mathcal{R}}_{\rho(f)}=m-1$,
then there exists some rotation $g'$ of $T^n$ satisfying
$$
\rho(g')=(\begin{array}{ccccccc}
\gamma_1&\gamma_2&\cdots&\gamma_m&0&\cdots&0
\end{array})^{\rm T},
$$
where $\gamma_1,\,\gamma_2,\,\cdots,\,\gamma_{m-1}$ are $m-1$ rationally independent irrational numbers,
but $\gamma_m=\dfrac{k}{d}$ is a rational number, where $d,\,k \in \mathbb{Z}$ and ${\rm gcd}\,(d, k)=1$,
such that $f$ and $g'$ are topologically conjugate.
\end{proposition}

\begin{proof}
According to the discussion in Remark \ref{rem:1},
assume that $\{\alpha_j\}_{j=1}^m$ is a maximal linearly independent subset of $\{\alpha_i\}_{i=1}^n$.
Naturally, $\alpha_{m+1}$ has a linear representation by $\alpha_1,\,\alpha_2,\,\cdots,\,\alpha_m$, i.e.,
$$
\alpha_{m+1}=r_1\alpha_1+r_2\alpha_2+\cdots+r_m\alpha_m,\quad r_1, r_2, \cdots, r_m \in \mathbb{Q}.
$$
So there exist some integers $k_1,\,k_2,\,\cdots,\,k_{m+1}$ satisfying $\gcd\,(k_1, k_2, \cdots, k_{m+1})=1$ such that
$$
k_1\alpha_1+k_2\alpha_2+\cdots+k_{m+1}\alpha_{m+1}=0.
$$
Lemma \ref{lem:6} indicates that there exists some matrix ${\bm A}_1 \in {\rm{GL}}_{m+1}(\mathbb{Z})$,
such that $(k_1,\,k_2,\,\cdots,\,k_{m+1})$ is just the $m+1$ row of ${\bm A}_1$.
Thus, we have
$$
\left(\begin{array}{c}
\alpha_{(1,\,1)}\\
\alpha_{(2,\,1)}\\
\vdots\\
\alpha_{(m,\,1)}\\
0\\
\alpha_{m+2}\\
\vdots\\
\alpha_n
\end{array}\right)=\left(\begin{array}{cc}
{\bm A}_1&0\\
0&{\bm I}_{n-m-1}
\end{array}\right)\left(\begin{array}{c}
\alpha_1\\
\alpha_2\\
\vdots\\
\alpha_m\\
\alpha_{m+1}\\
\alpha_{m+2}\\
\vdots\\
\alpha_n
\end{array}\right),
$$
where ${\bm I}_{n-m-1}$ is the identity matrix with order $n-m-1$.
Obviously, we can take a permutation matrix ${\bm T}_1$, such that
$$
{\bm T}_1\left(\begin{array}{c}
\alpha_{(1,\,1)}\\
\alpha_{(2,\,1)}\\
\vdots\\
\alpha_{(m,\,1)}\\
0\\
\alpha_{m+2}\\
\alpha_{m+3}\\
\vdots\\
\alpha_n
\end{array}\right)=
\left(\begin{array}{c}
\alpha_{(1,\,1)}\\
\alpha_{(2,\,1)}\\
\vdots\\
\alpha_{(m,\,1)}\\
\alpha_{m+2}\\
0\\
\alpha_{m+3}\\
\vdots\\
\alpha_n
\end{array}\right)={\bm T}_1\left(\begin{array}{cc}
{\bm A}_1&0\\
0&{\bm I}_{n-m-1}
\end{array}\right)\left(\begin{array}{c}
\alpha_1\\
\alpha_2\\
\vdots\\
\alpha_m\\
\alpha_{m+1}\\
\alpha_{m+2}\\
\alpha_{m+3}\\
\vdots\\
\alpha_n
\end{array}\right).
$$
Set ${\bm G}_1={\bm T}_1\left(\begin{array}{cc}
{\bm A}_1&0\\
0&{\bm I}_{n-m-1}
\end{array}\right)$.
One can see that ${\bm G}_1 \in {\rm{GL}}_n(\mathbb{Z})$.
So Proposition \ref{prop:1} implies that $\{\alpha_{(i, 1)}\}_{i=1}^m$ is a maximal linearly independent subset of $\{\alpha_{(1, 1)},\,\alpha_{(2, 1)}\,\cdots,\,\alpha_{(m, 1)},\,\alpha_{m+2},\,0,\,\alpha_{m+3},\,\cdots,\,\alpha_n\}$.
Then we proceed to do the same work by $n-m-1$ times. Finally, we can obtain
$$
\left(\begin{array}{c}
\alpha_{(1,\,n-m)}\\
\alpha_{(2,\,n-m)}\\
\vdots\\
\alpha_{(m,\,n-m)}\\
0\\
\vdots\\
0
\end{array}\right)={\bm G}_{n-m}{\bm G}_{n-m-1}\cdots{\bm G}_1\left(\begin{array}{c}
\alpha_1\\
\alpha_2\\
\vdots\\
\alpha_n
\end{array}\right),
$$
where ${\bm G}_i \in {\rm{GL}}_n(\mathbb{Z}),\,i=1,\,2,\,\cdots,\,n-m$.
Set
$$
\beta_i=\alpha_{(i,\,n-m)}\quad({\rm mod}\,\,\,\mathbb{Z}),\qquad i=1,\,2,\,\cdots,\,m.
$$
And take a rotation $g$ of $T^n$ satisfying
$$
\rho(g)=(\begin{array}{ccccccc}
\beta_1&\beta_2&\cdots&\beta_m&0&\cdots&0
\end{array})^{\rm T}.
$$
Thus, we know that $f$ and $g$ are topologically conjugate by introduction.

If $\widetilde{\mathcal{R}}_{\rho(f)}=\mathcal{R}_{\rho(f)}=m$, then according to the discussion in Remark \ref{rem:1},
it is obvious that $\beta_1,\,\beta_2,\,\cdots,\,\beta_m$ are $m$ rationally independent irrational numbers.
If $\widetilde{\mathcal{R}}_{\rho(f)}=\mathcal{R}_{\rho(f)}-1=m-1$,
then it follows from Definition \ref{rrr} and Proposition \ref{prop:1} that
$$
\widetilde{\mathcal{R}}(\beta_1, \beta_2, \cdots, \beta_m)=\widetilde{\mathcal{R}}(\beta_1, \beta_2, \cdots, \beta_m,\,0,\,\cdots,\,0)=\widetilde{\mathcal{R}}_{\rho(g)}=\widetilde{\mathcal{R}}_{\rho(f)}=m-1.
$$
Assume that $\{\beta_j\}_{j=1}^{m-1}$ is a maximal rationally independent subset of $\{\beta_i\}_{i=1}^m$.
Naturally, $\beta_{m}$ has a rational representation by $\beta_1,\,\beta_2,\,\cdots,\,\beta_{m-1}$, i.e.,
$$
\beta_m=\sum_{i=1}^{m-1}r_i\beta_i+r,\qquad r,\,\,\,r_i \in \mathbb{Q},\quad r \neq 0,
\quad i=1,\,2,\,\cdots,\,m-1.
$$
Then there exist some integers $k'_1,\,k'_2,\,\cdots,\,k'_m$ satisfying ${\rm gcd}\,(k'_1, k'_2, \cdots, k'_m)=1$ such that
$$
k'_1\beta_1+k'_2\beta_2+\cdots+k'_m\beta_m=\dfrac{d'}{d},
$$
where $d,\,d' \in \mathbb{Z}$, $d' \neq 0$ and ${\rm gcd}\,(d, d')=1$.
Lemma \ref{lem:6} implies that there exists some matrix ${\bm A}\in{\rm{GL}}_m(\mathbb{Z})$,
such that $(k'_1,\,k'_2,\,\cdots,\,k'_m)$ is just the $m$ row of ${\bm A}$.
Thus, we have
$$
\left(\begin{array}{c}
\beta'_1\\
\beta'_2\\
\vdots\\
\beta'_{m-1}\\
d'/d\\
0\\
\vdots\\
0
\end{array}\right)=\left(\begin{array}{cc}
{\bm A}&0\\
0&{\bm I}_{n-m}
\end{array}\right)\left(\begin{array}{c}
\beta_1\\
\beta_2\\
\vdots\\
\beta_{m-1}\\
\beta_m\\
0\\
\vdots\\
0
\end{array}\right),
$$
where ${\bm I}_{n-m}$ is the identity matrix with order $n-m$.
Take a rotation $g'$ of $T^n$ satisfying
$$
\rho(g')=(\begin{array}{ccccccc}
\gamma_1&\gamma_2&\cdots&\gamma_m&0&\cdots&0
\end{array})^{\rm T},
$$
where
$$
\gamma_i=\beta'_i\quad({\rm mod}\,\,\,\mathbb{Z}),
\qquad \gamma_m=\dfrac{k}{d}=\dfrac{d'}{d}\quad({\rm mod}\,\,\,\mathbb{Z}),
\qquad i=1,\,2,\,\cdots,\,m-1.
$$
Then according to Proposition \ref{prop:1}, we know
$$
\widetilde{\mathcal{R}}_{\rho(g)}=\widetilde{\mathcal{R}}_{\rho(g')}=\widetilde{\mathcal{R}}(\gamma_1, \gamma_2, \cdots, \gamma_m, 0,\cdots, 0)=m-1,
$$
i.e., $\gamma_1,\,\gamma_2,\,\cdots,\, \gamma_{m-1}$ are $m-1$ rationally independent irrational numbers,
and $\gamma_m=\dfrac{k}{d}$ is a rational number, where ${\rm gcd}\,(d, k)=1$.
Since $\left(\begin{array}{cc}
{\bm A}&0\\
0&{\bm I}_{n-m}
\end{array}\right) \in {\rm GL}_n(\mathbb{Z})$,
it follows from the discussion in introduction that $g$ and $g'$ are topologically conjugate.
Consequently, $f$ and $g'$ are topologically conjugate.
\end{proof}

\begin{proposition}\label{prop:3}
Assume that $f$ is a rotation of $T^n$ satisfying
$$
\rho(f)=(\begin{array}{cccc}
\alpha_1&\alpha_2&\cdots&\alpha_n
\end{array})^{\rm T},\qquad \alpha_i \in [0, 1),\,\,\,i=1, 2, \cdots, n,
$$
and $\widetilde{\mathcal{R}}_{\rho(f)}=m$.
Then $\overline{Orb_f(e)}$ is homeomorphic to the disjoint union of a group of $T^m$,
where $e$ is the identity element of $T^n$,
and
$$
\overline{Orb_f(e)}=\overline{\{e,\,f(e),\,f^2(e),\cdots,f^n(e),\cdots\}}
$$
is the orbit closure of $e$ under $f$.
\end{proposition}

\begin{proof}
Proposition \ref{prop:2} indicates that there exists some rotation $g$ of $T^n$ satisfying
$$
\rho(g)=(\begin{array}{cccccccc}
\beta_1&\beta_2&\cdots&\beta_m&\dfrac{k}{d}&0&\cdots&0
\end{array})^{\rm T},
$$
where $\beta_1,\,\beta_2,\,\cdots,\,\beta_m \in [0, 1)$ are $m$ rationally independent irrational numbers,
and $d,\,k \in \mathbb{Z},\,\,\,{\rm gcd}\,(d, k)=1$,
such that $f$ and $g$ are topologically conjugate.
According to lemma \ref{lem:1}, we know that there exists a topological conjugacy $h$ from $f$ to $g$ such that
$$
h \circ f=g \circ h,\qquad h(e)=e.
$$
Then one can see that
$$
h \circ f^n(e)=g^n \circ h(e),\qquad \forall\,n \in \mathbb{Z}_+.
$$
Thus, we have
$$
h(\overline{Orb_{f}(e)})=\overline{Orb_{g}(e)},
$$
i.e.,
$$
\overline{Orb_{f}(e)}\,\cong \overline{Orb_{g}(e)}.
$$
Notice that $g^d$ is also a rotation of $T^n$ satisfying
$$
\rho(g^d)=(\begin{array}{ccccccc}
\beta'_1&\beta'_2&\cdots&\beta'_m&0&\cdots&0
\end{array})^{\rm T},
$$
where
$$
\beta'_i=d\beta_i\quad({\rm mod}\,\,\,\mathbb{Z}),\qquad i=1,\,2,\,\cdots,\,m.
$$
Obviously, $\beta'_1,\,\beta'_2,\,\cdots,\,\beta'_m$ are $m$ rationally independent irrational numbers.
Set
$$
e_0=e,\quad e_1=g(e),\quad \cdots,\quad e_{d-1}=g^{d-1}(e).
$$
Then it follows from the form of $\rho(g^d)$ that we have
$$
\overline{Orb_g(e)}\,\cong\,\overline{Orb_{g^d}(e_0)}
\cup\overline{Orb_{g^d}(e_1)}\cup\cdots\cup\overline{Orb_{g^d}(e_{d-1})},
$$
$$
\overline{Orb_{g^d}(e_i)}\,\cong\,T^m,\qquad i=0,\,1,\,\cdots,\,d-1,
$$
and these $d$ orbit closures are mutually disjoint.
Therefore, $\overline{Orb_g(e)}$ is homeomorphic to the disjoint union of a group of $T^m$,
and then $\overline{Orb_f(e)}$ is homeomorphic to the disjoint union of a group of $T^m$.

In particular, when $d=1,\,k=0$, it is easy to see that $\overline{Orb_f(e)}$ is homeomorphic to (one) $T^m$.
So this case is obviously allowed.

\end{proof}

\section{Topologically conjugate classifications of the left actions on ${\rm SU}(2) \times T^n$}
In this section, for any positive integer $n$, we define the rotation vectors of the left actions in the set
$$
\mathscr{M}_{T_{{\rm SU}(2)\times T^n}}=\{L_g:\,{\rm SU}(2)\times T^n\rightarrow {\rm SU}(2)\times T^n;\,\, g \in T_{{\rm SU}(2)\times T^n}\},
$$
and utilize rotation vectors to give the topologically conjugate classifications of the left actions in $\mathscr{M}_{T_{{\rm SU}(2)\times T^n}}$,
and then give the topologically conjugate
classifications of all left actions on ${\rm SU}(2)\times T^n$.
Furthermore, we also study the relationship among their topological conjugacy, algebraic conjugacy and smooth conjugacy.

\subsection{Topologically conjugate classification theorem}

Set
$$
T_{{\rm SU}(2) \times T^n}=\left\{\left(\begin{array}{cc}
\lambda_1\cdots\lambda_n z&0\\
0&\bar{z}
\end{array}\right); z={{\textrm{e}}}^{2\pi{{\rm i}}\theta_0},
\lambda_i={{\textrm{e}}}^{2\pi{{\rm i}}\theta_i},\theta_0,\theta_1, \cdots,\theta_n\in [0, 1)\right\}.
$$
Obviously, $T_{{\rm SU}(2) \times T^n}$ is a maximal torus of ${\rm SU}(2) \times T^n$, and
$$
T_{{\rm SU}(2) \times T^n}\,\cong\,T_{{\rm SU}(2)} \times T^n\,\cong\,T^{n+1}.
$$
Define a map $\Phi:\,T_{{\rm SU}(2) \times T^n} \rightarrow T^{n+1}$ by
$$
\Phi:\,u=\left(\begin{array}{cc}
\lambda_1\lambda_2\cdots\lambda_n z&0\\
0&\bar{z}
\end{array}\right) \mapsto \left(\begin{array}{c}
z\\
\lambda_1\\
\vdots\\
\lambda_n
\end{array}\right),\qquad \forall\,u \in T_{{\rm SU}(2) \times T^n},
$$
It is easy to see that $\Phi$ is an isomorphism from the maximal torus $T_{{\rm SU}(2) \times T^n}$ to the normal $(n+1)$-dimensional torus $T^{n+1}$.
Then for any
$$
g=\left(\begin{array}{cc}
\lambda_1\lambda_2\cdots\lambda_n z&0\\
0&\bar{z}
\end{array}\right)\in T_{{\rm SU}(2) \times T^n},
\quad\hbox{i.e.},
\quad L_g \in \mathscr{M}_{T_{{\rm SU}(2) \times T^n}},
$$
set
$$
f=\Phi \circ L_g|_{T_{{\rm SU}(2) \times T^n}} \circ \Phi^{-1},
$$
where $z={{\textrm{e}}}^{2\pi{{\rm i}}\theta_0},
\,\lambda_i={{\textrm{e}}}^{2\pi{{\rm i}}\theta_i},\,\theta_i \in [0, 1),\,i=0,\,1,\,\cdots,\,n$.
One can see that $f:\,T^{n+1} \rightarrow T^{n+1}$ is a rotation (left action) of $T^{n+1}$ satisfying
$$
f:\,v \mapsto \left(\begin{array}{cccc}
{{\textrm{e}}}^{2\pi{{\rm i}}\theta_0}&&&\\
&{{\textrm{e}}}^{2\pi{{\rm i}}\theta_1}&&\\
&&\ddots&\\
&&&{{\textrm{e}}}^{2\pi{{\rm i}}\theta_n}
\end{array}\right)v,\qquad
\forall\,v\in T^{n+1}.
$$
This fact implies that $L_g|_{T_{{\rm SU}(2) \times T^n}}$ is topologically conjugate to some rotation $f$ of $T^{n+1}$.
Therefore, by Definition \ref{def:02}, we define the rotation vector of the left action $L_g$ under the representation $(T_{{\rm SU}(2) \times T^n}, \Phi)$ by
$$
\rho(L_g) \triangleq \rho(f)=(\begin{array}{cccc}
\theta_0&\theta_1&\cdots&\theta_n
\end{array})^{\rm T},
\qquad \theta_i \in [0, 1),\quad i=0,\,1,\,\cdots,\,n.
$$

Next, we give the topologically conjugate classification theorem.

\begin{theorem}\label{the:1}
For the left actions $L_g,\,L_{g'} \in \mathscr{M}_{T_{{\rm SU}(2) \times T^n}}$ with
$$
\rho(L_g)
=\left(\begin{array}{c}
\theta_0\\
\theta_1\\
\vdots\\
\theta_n
\end{array}\right),\qquad \rho(L_{g'})=\left(\begin{array}{c}
\theta'_0\\
\theta'_1\\
\vdots\\
\theta'_n
\end{array}\right),
$$
where $\theta_i,\,\theta'_i \in [0, 1),\,i=0, 1, \cdots, n$,
$L_g$ and $L_{g'}$ are topologically conjugate if and only if
$$
\left(\begin{array}{c}
\theta'_0\\
\theta'_1\\
\vdots\\
\theta'_n
\end{array}\right)=\left(\begin{array}{cc}
\pm\, 1& \bm{\ell}\\
\bm{0}&{\bm A}
\end{array}\right)\left(\begin{array}{c}
\theta_0\\
\theta_1\\
\vdots\\
\theta_n
\end{array}\right)\quad(\rm{mod}\,\,\,\mathbb{Z}^{n+1}),
$$
where $\bm{\ell}=(\begin{array}{cccc}\ell_1&\ell_2&\cdots&\ell_n\end{array})$ is a $1 \times n$ integer matrix, and ${\bm A} \in {{\rm GL}}_n(\mathbb{Z})$.
\end{theorem}

In fact, the equivalence relation in Theorem \ref{the:1} gives a complete topologically conjugate classification of the left actions in $\mathscr{M}_{T_{{\rm SU}(2) \times T^n}}$.
And then according to the discussion in introduction,
one can see that this equivalence relation gives a complete topologically conjugate classification of all left actions on ${\rm SU}(2) \times T^n$.

Furthermore, we give the relationship among the topological conjugacy, algebraic conjugacy and smooth conjugacy of the left actions on ${\rm SU}(2) \times T^n$ by following two propositions.

\begin{proposition}\label{prop:5}
There exist some left actions $L_g,\,L_{g'}$ on ${\rm SU}(2) \times T^n$ such that $L_g$ and $L_{g'}$ are topologically conjugate,
but not algebraically conjugate.
That means the topologically conjugate classification of the left actions on ${\rm SU}(2) \times T^n$ is not equivalent to their algebraically conjugate classification.
\end{proposition}

For the case $n=1$, the result has been proved in \cite{HP-2018}.
And when $n > 1$, it is not difficult to verify this proposition by the same way.

\begin{proposition}\label{prop:6}
For any left actions $L_g,\,L_{g'}$ on ${\rm SU}(2) \times T^n$,
$L_g$ and $L_{g'}$ are topologically conjugate if and only if $L_g$ and $L_{g'}$ are smooth conjugate.
That means the topologically conjugate classification of the left actions on Lie group ${\rm SU}(2) \times T^n$ is equivalent to their smooth conjugate classification.
\end{proposition}
The sufficiency of Proposition \ref{prop:6} is obviously true.
And we will prove the necessity of it at the end of this section.

\subsection{Proof of the sufficiency of Theorem \ref{the:1}}

In this part, we are going to construct topological conjugacies by the given relationship between rotation vectors to prove the sufficiency of Theorem \ref{the:1}.
For convenience, we denote the element $u=(z, {{\textrm{e}}}^{2\pi{{\rm i}}\varphi_1},{{\textrm{e}}}^{2\pi{{\rm i}}\varphi_2},\cdots,{{\textrm{e}}}^{2\pi{{\rm i}}\varphi_n}) \in {\rm SU}(2) \times T^n$ by
$$
u=(\begin{array}{ccccc}
z&\varphi_1&\varphi_2&\cdots&\varphi_n
\end{array})^{\rm T},
$$
where
$$
z=\left(\begin{array}{cc}
z_1&-\bar{z_2}\\
z_2&\bar{z_1}
\end{array}\right) \in {\rm SU}(2),
\,\, z_1,\,z_2 \in \mathbb{C},\,\,\,|z_1|^2+|z_2|^2=1,
\,\,\varphi_1,\,\varphi_2,\,\cdots,\,\varphi_n \in \mathbb{R}.
$$

\begin{proof}
We know that the left actions $L_g,\,L_{g'} \in \mathscr{M}_{T_{{\rm SU}(2) \times T^n}}$ satisfy
$$
\left(\begin{array}{c}
\theta'_0\\
\theta'_1\\
\vdots\\
\theta'_n
\end{array}\right)=\rho(L_{g'})=\left(\begin{array}{cc}
\pm\, 1& \bm{\ell}\\
\bm{0}&{\bm A}
\end{array}\right)\rho(L_g)=\left(\begin{array}{cc}
\pm\, 1& \bm{\ell}\\
\bm{0}&{\bm A}
\end{array}\right)\left(\begin{array}{c}
\theta_0\\
\theta_1\\
\vdots\\
\theta_n
\end{array}\right)\,\,\,(\rm{mod}\,\,\mathbb{Z}^{n+1}),
$$
where $\theta_i,\,\theta'_i \in [0, 1),\,i=0, 1, \cdots, n$, ${\bm A} \in {{\rm GL}}_n(\mathbb{Z})$,
and $\bm{\ell}=(\begin{array}{cccc}\ell_1&\ell_2&\cdots&\ell_n\end{array})$ is a $1 \times n$ integer matrix.
Thus, there are two different cases.

(1) Assume that
$$
\left(\begin{array}{c}
\theta'_0\\
\theta'_1\\
\vdots\\
\theta'_n
\end{array}\right)=\left(\begin{array}{cc}
1& \bm{\ell}\\
\bm{0}&{\bm A}
\end{array}\right)\left(\begin{array}{c}
\theta_0\\
\theta_1\\
\vdots\\
\theta_n
\end{array}\right)\quad(\rm{mod}\,\,\,\mathbb{Z}^{n+1}).
$$
Then for any
$$
u=\left(\begin{array}{ccccc}
\left(\begin{array}{cc}
z_1&-\bar{z_2}\\
z_2&\bar{z_1}
\end{array}\right)&\varphi_1&\varphi_2&\cdots&\varphi_n
\end{array}\right)^{\rm T} \in {\rm SU}(2) \times T^n,
$$
define the maps $h_1,\,h_1':\,{\rm SU}(2) \times T^n \rightarrow {\rm SU}(2) \times T^n$ by
\begin{gather*}
h_1:\,u=\left(\begin{array}{c}
\left(\begin{array}{cc}
z_1&-\bar{z_2}\\
z_2&\bar{z_1}
\end{array}\right)\\
{\bm \varphi}
\end{array}\right) \mapsto \left(\begin{array}{c}
{\left(\begin{array}{cc}
{{\textrm{e}}}^{2\pi{{\rm i}}{\bm \ell}{\bm \varphi}}&0\\
0&{{\textrm{e}}}^{-2\pi{{\rm i}}{\bm \ell}{\bm \varphi}}
\end{array}\right)}{\left(\begin{array}{cc}
z_1&-\bar{z_2}\\
z_2&\bar{z_1}
\end{array}\right)}\\
{\bm A}{\bm \varphi}
\end{array}\right),
\\
h_1':\,u=\left(\begin{array}{c}
\left(\begin{array}{cc}
z_1&-\bar{z_2}\\
z_2&\bar{z_1}
\end{array}\right)\\
{\bm \varphi}
\end{array}\right) \mapsto \left(\begin{array}{c}
{\left(\begin{array}{cc}
{{\textrm{e}}}^{-2\pi{{\rm i}}{\bm \ell}{\bm A}^{-1}{\bm \varphi}}&0\\
0&{{\textrm{e}}}^{2\pi{{\rm i}}{\bm \ell}{\bm A}^{-1}{\bm \varphi}}
\end{array}\right)}{\left(\begin{array}{cc}
z_1&-\bar{z_2}\\
z_2&\bar{z_1}
\end{array}\right)}\\
{\bm A}^{-1}{\bm \varphi}
\end{array}\right).
\end{gather*}
where ${\bm \varphi}=(\begin{array}{cccc}
\varphi_1&\varphi_2&\cdots&\varphi_n
\end{array})^{\rm T}$.
One can see that $h_1$ and $h_1'$ are both continuous, and $h_1' \circ h_1={\rm id}_{{\rm SU}(2) \times T^n}$,
where ${\rm id}_{{\rm SU}(2) \times T^n}$ is the identity map of ${\rm SU}(2) \times T^n$.
Then $h_1$ is a self-homeomorphism of ${\rm SU}(2) \times T^n$.
It follows from a simple verification that we have
$$
h_1 \circ L_g(v)=L_{g'} \circ h_1(v),\qquad \forall\,v  \in {\rm SU}(2) \times T^n.
$$
Therefore, $L_g$ and $L_{g'}$ are topologically conjugate.

(2) Assume that
$$
\left(\begin{array}{c}
\theta'_0\\
\theta'_1\\
\vdots\\
\theta'_n
\end{array}\right)=\left(\begin{array}{cc}
-1& \bm{\ell}\\
\bm{0}&{\bm A}
\end{array}\right)\left(\begin{array}{c}
\theta_0\\
\theta_1\\
\vdots\\
\theta_n
\end{array}\right)\quad(\rm{mod}\,\,\,\mathbb{Z}^{n+1}).
$$
Then for any
$$
u=\left(\begin{array}{ccccc}
\left(\begin{array}{cc}
z_1&-\bar{z_2}\\
z_2&\bar{z_1}
\end{array}\right)&\varphi_1&\varphi_2&\cdots&\varphi_n
\end{array}\right)^{\rm T} \in {\rm SU}(2) \times T^n,
$$
define the maps $h_2,\,h_2':\,{\rm SU}(2) \times T^n \rightarrow {\rm SU}(2) \times T^n$ by
$$
h_2:\,u=\left(\begin{array}{c}
\left(\begin{array}{cc}
z_1&-\bar{z_2}\\
z_2&\bar{z_1}
\end{array}\right)\\
{\bm \varphi}
\end{array}\right) \mapsto \left(\begin{array}{c}
{\left(\begin{array}{cc}
{{\textrm{e}}}^{2\pi{{\rm i}}{\bm \ell}{\bm \varphi}}&0\\
0&{{\textrm{e}}}^{-2\pi{{\rm i}}{\bm \ell}{\bm \varphi}}
\end{array}\right)}{\left(\begin{array}{cc}
\bar{z_1}&-z_2\\
\bar{z_2}&z_1
\end{array}\right)}\\
{\bm A}{\bm \varphi}
\end{array}\right),
$$
$$
h_2':\,u=\left(\begin{array}{c}
\left(\begin{array}{cc}
z_1&-\bar{z_2}\\
z_2&\bar{z_1}
\end{array}\right)\\
{\bm \varphi}
\end{array}\right) \mapsto \left(\begin{array}{c}
{\left(\begin{array}{cc}
{{\textrm{e}}}^{2\pi{{\rm i}}{\bm \ell}{\bm A}^{-1}{\bm \varphi}}&0\\
0&{{\textrm{e}}}^{-2\pi{{\rm i}}{\bm \ell}{\bm A}^{-1}{\bm \varphi}}
\end{array}\right)}{\left(\begin{array}{cc}
\bar{z_1}&-z_2\\
\bar{z_2}&z_1
\end{array}\right)}\\
{\bm A}^{-1}{\bm \varphi}
\end{array}\right),
$$
respectively, where ${\bm \varphi}=(\begin{array}{cccc}
\varphi_1&\varphi_2&\cdots&\varphi_n
\end{array})^{\rm T}$.
Similar to Case (1), we know that $h_2$ is a self-homeomorphism of ${\rm SU}(2) \times T^n$ satisfying
$$
h_2 \circ L_g(v)=L_{g'} \circ h_2(v)\qquad \forall\,v \in {\rm SU}(2) \times T^n.
$$
So $L_g$ and $L_{g'}$ are topologically conjugate.
\end{proof}

\subsection{Proof of the necessity of Theorem \ref{the:1}}

In this part, we utilize the relevant concepts and results in Section $2$ to give the proof of the necessity of Theorem \ref{the:1}.

\begin{proof}
Assume that the left actions $L_g,\,L_{g}' \in \mathscr{M}_{T_{{\rm SU}(2) \times T^n}}$ are topologically conjugate.
Then Lemma \ref{lem:1} indicates that there exists a topological conjugacy from $L_g$ to $L_{g'}$ such that
$$
h \circ L_g=L_{g'} \circ h,\qquad h(e)=e,
$$
where $e$ is the identity element of ${\rm SU}(2) \times T^n$.
And we claim that the topological conjugacies in the present paper always preserve the identity elements.

Investigate the following diagram.
$$
\xymatrixcolsep{3pc}
\xymatrix{
{\,\,T^{n+1}\,\,} \ar[d]_-{f_1} \ar[r]^-{i}
& {\,\,{\rm SU}(2) \times T^n\,\,} \ar[d]_-{L_g} \ar[r]^-{h}
& {\,\,{\rm SU}(2) \times T^n\,\,} \ar[d]_-{L_{g'}} \ar[r]^-{\pi}
& {\,\,T^n\,\,} \ar[d]_-{f_2}\\
{\,\,T^{n+1}\,\,} \ar[r]_-{i}
& {\,\,{\rm SU}(2) \times T^n\,\,} \ar[r]_-{h}
& {\,\,{\rm SU}(2) \times T^n\,\,} \ar[r]_-{\pi}
& {\,\,T^n\,\,}.}
$$
In the diagram, $i$ is the natural inclusion preserving identity elements, $\pi$ is the projection satisfying
$$
\pi:\,(u, t) \mapsto t,\qquad \forall\,u \in {\rm SU}(2),\,t \in T^n,
$$
$f_1,\,f_2$ are the rotations of $T^{n+1}$ and $T^n$, respectively,
$L_g,\,L_{g'}\in\mathscr{M}_{T_{{\rm SU}(2) \times T^n}}$ are two left actions on ${\rm SU}(2) \times T^n$, $h$ is the topological conjugacy from $L_g$ to $L_{g'}$,
and $f_1,\,f_2,\,L_g,\,L_{g'}$ satisfy
$$
\rho(f_1)=\rho(L_g)=\left(\begin{array}{c}
\theta_0\\
\theta_1\\
\vdots\\
\theta_n
\end{array}\right),
\qquad \rho(L_{g'})=\left(\begin{array}{c}
\theta'_0\\
\theta'_1\\
\vdots\\
\theta'_n
\end{array}\right),\qquad \rho(f_2)=\left(\begin{array}{c}
\theta'_1\\
\theta'_2\\
\vdots\\
\theta'_n
\end{array}\right).
$$
It is easy to know that
$$
i \circ f_1=L_g \circ i,\qquad \pi \circ L_{g'}=f_2 \circ \pi.
$$
Thus, the above diagram is commutative,
and then we have
$$
f \circ f_1=f_2 \circ f,
$$
where $f=\pi \circ h \circ i$ is a continuous surjection from $T^{n+1}$ to $T^n$.
It is well-known that
$$
\pi_1({\rm SU}(2) \times T^n)\,\cong\,\pi_1(T^n)\,\cong \,\mathbb{Z}^n.
$$
Then for any
$$
{\bm l}=\left(\begin{array}{c}
l_0\\
l_1\\
\vdots\\
l_n
\end{array}\right) \in \pi_1(T^{n+1})\,\cong\,\mathbb{Z}^{n+1},
\qquad {\bm l}'=\left(\begin{array}{c}
l'_1\\
l'_2\\
\vdots\\
l'_n
\end{array}\right) \in \pi_1({\rm SU}(2) \times T^n) \cong \mathbb{Z}^n,
$$
the homomorphisms $i_*,\,h_*,\,\pi_*$ between fundamental groups induced by $i,\,h,\,\pi$, respectively, satisfy
$$
i_*({\bm l})={\bm A}_1{\bm l}=\left(\begin{array}{c}
l_1\\
l_2\\
\vdots\\
l_n
\end{array}\right),
\qquad h_*({\bm l}')={\bm A}{\bm l}',
\qquad \pi_*({\bm l}')={\bm A}_2{\bm l}'=\left(\begin{array}{c}
l'_1\\
l'_2\\
\vdots\\
l'_n
\end{array}\right),
$$
where ${\bm A}_1=\left(\begin{array}{cc}
{\bm 0}&{\bm I}_n
\end{array}\right)$ is an $n \times (n+1)$ integer matrix, ${\bm A} \in {\rm GL}_n(\mathbb{Z})$ is the matrix form of the isomorphism $h_*$,
${\bm A}_2={\bm I}_n$ is the identity matrix with order $n$.
Thus, the homomorphism $f_*:\,\pi_1(T^{n+1})  \rightarrow \pi_1(T^n)$ induced by $f=\pi \circ h \circ i$ satisfies
$$
f_*({\bm l})=(\pi \circ h \circ i)_*({\bm l})
=\pi_* \circ h_* \circ i_*({\bm l})={\bm A}_2{\bm A}{\bm A}_1{\bm l},
\quad \forall\,{\bm l} \in \pi_1(T^{n+1}).
$$
One can see that ${\bm A}_2{\bm A}{\bm A}_1$ is an $n \times (n+1)$ integer matrix satisfying
$$
{\bm A}_2{\bm A}{\bm A}_1=(\begin{array}{cc}
{\bm 0}&{\bm A}
\end{array}).
$$
Finally, it follows from Lemma \ref{lem:2} that we get
$$
\rho(f_2)=(\begin{array}{cc}
{\bm 0}&{\bm A}
\end{array})\,\rho(f_1)\quad({\rm mod}\,\,\,\mathbb{Z}^n),
$$
i.e.,
$$
\left(\begin{array}{c}
\theta'_1\\
\theta'_2\\
\vdots\\
\theta'_n
\end{array}\right)=(\begin{array}{cc}
{\bm 0}&{\bm A}
\end{array})\left(\begin{array}{c}
\theta_0\\
\theta_1\\
\vdots\\
\theta_n
\end{array}\right)\quad(\rm{mod}\,\,\,\mathbb{Z}^n).
$$
Therefore, it suffices to prove
$$
\theta'_0=\pm\,\theta_0+\sum_{i=1}^n{\ell_i\theta_i}\quad({\rm{mod}}\,\,\,\mathbb{Z}),
\qquad \ell_i \in \mathbb{Z},\,\,\, i=1,\,2,\,\cdots,\,n.
$$

Take $L_g\ \in \mathscr{M}_{T_{{\rm SU}(2) \times T^n}}$ satisfying
$$
\rho(L_g)=(\begin{array}{cccc}
\theta_0&\theta_1&\cdots&\theta_n
\end{array})^{\rm T},\qquad \theta_i \in [0, 1),\,\,\, i=0,\,1,\,\cdots,\,n.
$$
And define
$$
\widetilde{\mathcal{R}}_{\rho(L_g)}\triangleq\widetilde{\mathcal{R}}(\theta_0, \theta_1,\cdots, \theta_n),
\qquad\widetilde{\mathcal{R}}'_{\rho(L_g)}\triangleq\widetilde{\mathcal{R}}(\theta_1, \theta_2,\cdots, \theta_n).
$$
Assume that
$$
\widetilde{\mathcal{R}}'_{\rho(L_g)}=\widetilde{\mathcal{R}}(\theta_1, \theta_2,\cdots, \theta_n)=m.
$$
Then Proposition \ref{prop:2} indicates that there exist $m$ rationally independent irrational numbers $\beta_1,\,\beta_2,\,\cdots,\,\beta_m \in [0, 1)$ and a rational number $\dfrac{k}{d} \in [0, 1)$,
where $d,\,k \in \mathbb{Z}$ and ${\rm gcd}\,(d, k)=1$, such that
$$
\left(\begin{array}{c}
\beta_1\\
\vdots\\
\beta_m\\
k/d\\
0\\
\vdots\\
0
\end{array}\right)={\bm A}_g\left(\begin{array}{c}
\theta_1\\
\theta_2\\
\vdots\\
\theta_n\\
\end{array}\right)\quad({\rm{mod}}\,\,\,\mathbb{Z}^n),
\qquad {\bm A}_g \in {\rm GL}_n(\mathbb Z),
$$
Thus, we have
$$
\left(\begin{array}{c}
\theta_0\\
\beta_1\\
\vdots\\
\beta_m\\
k/d\\
0\\
\vdots\\
0
\end{array}\right)=\left(\begin{array}{cc}
1&{\bm 0}\\
{\bm 0}&{\bm A}_g
\end{array}\right)\left(\begin{array}{c}
\theta_0\\
\theta_1\\
\vdots\\
\theta_n\\
\end{array}\right)\quad({\rm{mod}}\,\,\,\mathbb{Z}^{n+1}),
\qquad {\bm A}_g \in {\rm GL}_n(\mathbb Z).
$$
Take $L_t \in \mathscr{M}_{T_{{\rm SU}(2) \times T^n}}$ satisfying
$$
\rho(L_t)=(\begin{array}{cccccccc}
\theta_0&\beta_1&\cdots&\beta_m&k/d&0&\cdots&0
\end{array})^{\rm T}.
$$
Then the sufficiency of Theorem \ref{the:1} implies that $L_g$ and $L_t$ are topologically conjugate.
So we can identify $L_g$ with $L_t$ up to topologically conjugate equivalence.
In this way, we say that $\rho(L_g)$ possesses a normal form
$$
\rho(L_t)=(\begin{array}{cccccccc}
\theta_0&\beta_1&\cdots&\beta_m&k/d&0&\cdots&0
\end{array})^{\rm T}.
$$
If $L_g$ and $L_{g'}$ are topologically conjugate,
and $\rho(L_t),\,\rho(L_{t'})$ are the normal forms of $\rho(L_g)$ and $\rho(L_{g'})$, respectively,
then $L_t$ and $L_{t'}$ are topologically conjugate.
Assume that $\rho(L_t)$ and $\rho(L_{t'})$ satisfy
$$
\rho(L_{t'})=\left(\begin{array}{cc}
\pm\,1&{\bm \ell}'\\
{\bm {0}}&{\bm A}'
\end{array}\right)\rho(L_t)\quad({\rm{mod}}\,\,\,\mathbb{Z}^{n+1}).
$$
Then by the above discussion, it is easy to prove that
$$
\rho(L_{g'})=\left(\begin{array}{cc}
\pm\,1&{\bm \ell}\\
{\bm {0}}&{\bm A}
\end{array}\right)\rho(L_g)\quad({\rm{mod}}\,\,\,\mathbb{Z}^{n+1}),
$$
where ${\bm \ell},\,{\bm \ell}'$ are $1 \times n$ integer matrices,
and ${\bm A},\,{\bm A}' \in {{\rm GL}}_n(\mathbb{Z})$.
Consequently, in the following proof, we just need to study the left actions satisfying
$$
\rho(L_g)=(\begin{array}{cccccccc}
\theta_0&\beta_1&\cdots&\beta_m&k/d&0&\cdots&0
\end{array})^{\rm T}
$$
when $\widetilde{\mathcal{R}}'_{\rho(L_g)}=m$.

\begin{remark}\label{rem:3}
Notice that we fix $\theta_0$ in the above discussion about the normal forms of rotation vectors.
Since the previous proof has illuminated that if $L_g$ and $L_{g'}$ satisfying
$$
\rho(L_g)
=(\begin{array}{cccc}
\theta_0&\theta_1&\cdots&\theta_n
\end{array})^{\rm T},\qquad \rho(L_{g'})=(\begin{array}{cccc}
\theta'_0&\theta'_1&\cdots&\theta'_n
\end{array})^{\rm T}
$$
are topologically conjugate,
we have
$$
\left(\begin{array}{c}
\theta'_1\\
\theta'_2\\
\vdots\\
\theta'_n
\end{array}\right)=(\begin{array}{cc}
{\bm 0}&{\bm A}
\end{array})\left(\begin{array}{c}
\theta_0\\
\theta_1\\
\vdots\\
\theta_n
\end{array}\right)\quad({\rm mod}\,\,\,\mathbb{Z}^n),
\qquad {\bm A} \in {\rm GL}_n(\mathbb{Z}),
$$
then reviewing the reduction process of rotation vectors in the proof of Proposition \ref{prop:2},
one can see that it is necessary to fix $\theta_0$.
\end{remark}

Assume that $L_g,\,L_{g'} \in \mathscr{M}_{T_{{\rm SU}(2) \times T^n}}$ satisfy
$$
\rho(L_g)=(\begin{array}{cccc}
\theta_0&\theta_1&\cdots&\theta_n
\end{array})^{\rm T},\qquad \rho(L_{g'})=(\begin{array}{cccc}
\theta'_0&\theta'_1&\cdots&\theta'_n
\end{array})^{\rm T},
$$
and $h$ is a topological conjugacy from $L_g$ to $L_{g'}$ satisfying $h(e)=e$,
where $e$ is the identify element of ${\rm SU}(2) \times T^n$.
Then
$$
h \circ L_g^n(e)=L_{g'}^n \circ h(e),
$$
As $n \rightarrow \infty$, we have
$$
h(\overline{Orb_{L_g}(e)})=\overline{Orb_{L_{g'}}(e)},
\quad\hbox{i.e.,}\quad\overline{Orb_{L_g}(e)}\cong\overline{Orb_{L_{g'}}(e)}.
$$
According to the previous analysis, we know that the left actions $L_g|_{T_{{\rm SU}(2) \times T^n}}$ and $L_{g'}|_{T_{{\rm SU}(2) \times T^n}}$ are topologically conjugate to the rotations $f$ and $f'$ of $T^{n+1}$ with
$$
\rho(f)=(\begin{array}{cccc}
\theta_0&\theta_1&\cdots&\theta_n
\end{array})^{\rm T},\qquad \rho(f')=(\begin{array}{cccc}
\theta'_0&\theta'_1&\cdots&\theta'_n
\end{array})^{\rm T},
$$
respectively.
Similarly, there also exist some topological conjugacies preserving the identify elements between them.
So we have
$$
\overline{Orb_{L_g|_{T_{{\rm SU}(2) \times T^n}}}(e)}\,\cong\,\overline{Orb_f(e_{n+1})},
\qquad \overline{Orb_{L_{g'}|_{T_{{\rm SU}(2) \times T^n}}}(e)}\,\cong\,\overline{Orb_{f'}(e_{n+1})},
$$
where $e_{n+1}$ is the identify element of $T^{n+1}$.
Since $L_g,\,L_{g'} \in \mathscr{M}_{T_{{\rm SU}(2) \times T^n}}$, and $e \in T_{{\rm SU}(2) \times T^n}$, then
$$
\overline{Orb_{L_g}(e)}=\overline{Orb_{L_g|_{T_{{\rm SU}(2) \times T^n}}}(e)},
\qquad \overline{Orb_{L_{g'}}(e)}=\overline{Orb_{L_{g'}|_{T_{{\rm SU}(2) \times T^n}}}(e)}.
$$
Thus, we obtain
$$
\overline{Orb_{L_g|_{T_{{\rm SU}(2) \times T^n}}}(e)}\,\cong\,\overline{Orb_{L_{g'}|_{T_{{\rm SU}(2) \times T^n}}}(e)},
\quad\hbox{i.e.,}\quad \overline{Orb_f(e_{n+1})}\,\cong\,\overline{Orb_{f'}(e_{n+1})}.
$$
It follows from Proposition \ref{prop:3} that $\overline{Orb_f(e_{n+1})}$ and $\overline{Orb_{f'}(e_{n+1})}$ are both homeomorphic to the disjoint union of a group of $T^{m'}$,
i.e., $\overline{Orb_{L_g}(e)}$ and $\overline{Orb_{L_{g'}}(e)}$ are both homeomorphic to the disjoint union of a group of $T^{m'}$,
and the reductive ranks
$$
\widetilde{\mathcal{R}}_{\rho(L_g)}=\widetilde{\mathcal{R}}_{\rho(f)},
\qquad \widetilde{\mathcal{R}}_{\rho(L_{g'})}=\widetilde{\mathcal{R}}_{\rho(f')}
$$
are just equal to the dimensions of $\overline{Orb_{L_g}(e)}$ and $\overline{Orb_{L_{g'}}(e)}$, respectively.
Thus, we get
$$
\widetilde{\mathcal{R}}_{\rho(L_g)}=\widetilde{\mathcal{R}}_{\rho(L_{g'})}=m'.
$$
This fact means $\widetilde{\mathcal{R}}$ is a topologically conjugate invariant.
Moreover, we denote the number of the mutually disjoint $T^{m'}$ in $\overline{Orb_{L_g}(e)}$ by $\mathcal{N}_{\rho(L_g)}$,
and assume that $\mathcal{N}_{\rho(L_g)}=l$.
Then it is obvious that
$$
\mathcal{N}_{\rho(L_g)}=\mathcal{N}_{\rho(L_{g'})}=l.
$$
That means $\mathcal{N}$ is also a topologically conjugate invariant.

On the other hand, let $f_0,\,f'_0$ be two rotations of $T^n$ with
$$
\rho(f_0)=(\begin{array}{cccc}
\theta_1&\theta_2&\cdots&\theta_n
\end{array})^{\rm T},\qquad \rho(f'_0)=(\begin{array}{cccc}
\theta'_1&\theta'_2&\cdots&\theta'_n
\end{array})^{\rm T}.
$$
The previous analysis indicates that
$$
\left(\begin{array}{c}
\theta'_1\\
\theta'_2\\
\vdots\\
\theta'_n
\end{array}\right)=(\begin{array}{cc}
{\bm 0}&{\bm A}
\end{array})\left(\begin{array}{c}
\theta_0\\
\theta_1\\
\vdots\\
\theta_n
\end{array}\right)\quad({\rm mod}\,\,\,\mathbb{Z}^n),
\qquad {\bm A} \in {\rm GL}_n(\mathbb{Z}),
$$
i.e.,
$$
\left(\begin{array}{c}
\theta'_1\\
\theta'_2\\
\vdots\\
\theta'_n
\end{array}\right)={\bm A}\left(\begin{array}{c}
\theta_1\\
\theta_2\\
\vdots\\
\theta_n
\end{array}\right)\quad({\rm mod}\,\,\,\mathbb{Z}^n),
\qquad {\bm A} \in {\rm GL}_n(\mathbb{Z}).
$$
It follows from the discussion in introduction that $f_0$ and $f'_0$ are topologically conjugate.
Take the topological conjugacy $h'$ from $f_0$ to $f'_0$ satisfying $h'(e_n)=e_n$, where $e_n$ is the identity element of $T^n$.
By the same way, one can see that $\overline{Orb_{f_0}(e_n)}$ and $\overline{Orb_{f'_0}(e_n)}$ are both homeomorphic to the disjoint union of a group of $T^m$,
and the reductive ranks $\widetilde{\mathcal{R}}_{\rho(f_0)},\,\widetilde{\mathcal{R}}_{\rho(f'_0)}$ are just equal to the dimensions of $\overline{Orb_{f_0}(e_n)}$ and $\overline{Orb_{f'_0}(e_n)}$, respectively.
In fact, one can see that
$$
\widetilde{\mathcal{R}}'_{\rho(L_g)}=\widetilde{\mathcal{R}}_{\rho(f_0)},
\qquad \widetilde{\mathcal{R}}'_{\rho(L_{g'})}=\widetilde{\mathcal{R}}_{\rho(f'_0)}
$$
Thus, we have
$$
\widetilde{\mathcal{R}}'_{\rho(L_g)}=\widetilde{\mathcal{R}}'_{\rho(L_{g'})}=m.
$$
This fact implies that $\widetilde{\mathcal{R}}'$ is a topologically conjugate invariant.
Moreover, we denote the number of the mutually disjoint $T^m$ in $\overline{Orb_{f_0}(e_n)}$ by $\mathcal{N}'_{\rho(L_g)}$,
and assume that $\mathcal{N}'_{\rho(L_g)}=d$.
Then it is obvious that
$$
\mathcal{N}'_{\rho(L_g)}=\mathcal{N}'_{\rho(L_{g'})}=d.
$$
That means $\mathcal{N}'$ is also a topologically conjugate invariant.
To sum up, for the left actions in $\mathscr{M}_{T_{{\rm SU}(2) \times T^n}}$,
we define $4$ topologically conjugate invariants about them,
i.e., $\widetilde{\mathcal{R}},\,\widetilde{\mathcal{R}}',\,\mathcal{N},\,\mathcal{N}'$.

Assume that $L_g,\,L_{g'} \in \mathscr{M}_{T_{{\rm SU}(2) \times T^n}}$ are topologically conjugate.
Then combining the proof of Proposition \ref{prop:3} with the discussions about normal forms and topologically conjugate invariants,
we set
$$
\rho(L_g)=\left(\begin{array}{c}
\theta_0\\
\beta_1\\
\vdots\\
\beta_m\\
k/d\\
0\\
\vdots\\
0
\end{array}\right),
\qquad\qquad \rho(L_{g'})=\left(\begin{array}{c}
\theta'_0\\
\beta'_1\\
\vdots\\
\beta'_m\\
k'/d\\
0\\
\vdots\\
0
\end{array}\right),
$$
where $\{\beta_i\}_{i=1}^m$ and $\{\beta'_j\}_{j=1}^m$ are both $m$ rational independent irrational numbers,
$d,\,k,\,k' \in \mathbb{Z}$ and ${\rm gcd}\,(d, k)={\rm gcd}\,(d, k')=1$.
One can see that
$$
\widetilde{\mathcal{R}}'_{\rho(L_g)}=\widetilde{\mathcal{R}}'_{\rho(L_{g'})}=m,
\qquad \mathcal{N}'_{\rho(L_g)}=\mathcal{N}'_{\rho(L_{g'})}=d.
$$
Then for the left action $L_g$,
it is no hard to see that there exist 3 possible cases about the topologically conjugate invariants $\widetilde{\mathcal{R}}_{\rho(L_g)},\,\mathcal{N}_{\rho(L_g)},\,\widetilde{\mathcal{R}}'_{\rho(L_g)},\,\mathcal{N}'_{\rho(L_g)}$ as follows.

{\bf Case 1.} $\theta_0\in [0, 1)$ is a rational number. Then
$$
\widetilde{\mathcal{R}}_{\rho(L_g)}=m,
\qquad \mathcal{N}_{\rho(L_g)}=l,
\qquad \widetilde{\mathcal{R}}'_{\rho(L_g)}=m,
\qquad \mathcal{N}'_{\rho(L_g)}=d,
\qquad d \leq l \in \mathbb{Z}.
$$

{\bf Case 2.} $\theta_0$ is an irrational number, and $\theta_0,\,\beta_1,\,\beta_2,\,\cdots,\,\beta_m$ are rational independent. Then
$$
\widetilde{\mathcal{R}}_{\rho(L_g)}=m+1,
\qquad \mathcal{N}_{\rho(L_g)}=d,
\qquad \widetilde{\mathcal{R}}'_{\rho(L_g)}=m,
\qquad \mathcal{N}'_{\rho(L_g)}=d.
$$

{\bf Case 3.} $\theta_0$ is an irrational number, and $\theta_0,\,\beta_1,\,\beta_2,\,\cdots,\,\beta_m$ are rational dependent. Then
$$
\widetilde{\mathcal{R}}_{\rho(L_g)}=m,
\qquad \mathcal{N}_{\rho(L_g)}=l,
\qquad \widetilde{\mathcal{R}}'_{\rho(L_g)}=m,
\qquad \mathcal{N}'_{\rho(L_g)}=d,
\qquad d \leq l \in \mathbb{Z}.
$$

Obviously, for the left action $L_{g'}$, the possible cases are the same as {\bf Case 1}, {\bf Case 2} and {\bf Case 3}.
Therefore, we divide the following proof into $4$ part by the topologically conjugate invariance of $\widetilde{\mathcal{R}},\,\mathcal{N},\,\widetilde{\mathcal{R}}',\,\mathcal{N}'$,
i.e.,
$$
\widetilde{\mathcal{R}}_{\rho(L_g)}=\widetilde{\mathcal{R}}_{\rho(L_{g'})},
\quad \mathcal{N}_{\rho(L_g)}=\mathcal{N}_{\rho(L_{g'})},
\quad \widetilde{\mathcal{R}}'_{\rho(L_g)}=\widetilde{\mathcal{R}}'_{\rho(L_{g'})},
\quad \mathcal{N}'_{\rho(L_g)}=\mathcal{N}'_{\rho(L_{g'})}.
$$

${\bm {(1)}}$ Assume that both of $\rho(L_g)$ and $\rho(L_g)$ satisfy the conditions in {\bf Case 1}.

Let $\theta_0,\,\theta'_0 \in [0, 1)$ be rational numbers, and
$$
\dfrac{q}{p}=d\theta_0\,\,\,({\rm mod}\,\,\,\mathbb{Z}),
\qquad \dfrac{q'}{p'}=d\theta'_0\,\,\,({\rm mod}\,\,\,\mathbb{Z}),
\qquad {\rm gcd}\,(p, q)={\rm gcd}\,(p', q')=1.
$$
And take the topological conjugacy $h$ from $L_g$ to $L_{g'}$ preserving the identify element $e$ of ${\rm SU}(2) \times T^n$.
Then we have
$$
h \circ L_g^d=L_{g'}^d \circ h,\qquad h(\overline{Orb_{L_g^d}(e)})=\overline{Orb_{L_{g'}^d}(e)},
$$
and $L_g^d,\,L_{g'}^d \in \mathscr{M}_{T_{{\rm SU}(2) \times T^n}}$ with
$$
\rho(L_g^d)=\left(\begin{array}{c}
q/p\\
d\beta_1\\
\vdots\\
d\beta_m\\
k\\
0\\
\vdots\\
0
\end{array}\right)\quad({\rm mod}\,\,\,\mathbb{Z}^{n+1}),
\qquad \rho(L_{g'}^d)=\left(\begin{array}{c}
q'/p'\\
d\beta'_1\\
\vdots\\
d\beta'_m\\
k'\\
0\\
\vdots\\
0
\end{array}\right)\quad({\rm mod}\,\,\,\mathbb{Z}^{n+1}),
$$
where $\{d\beta_i\}_{i=1}^m$ and $\{d\beta'_i\}_{i=1}^m$ are both $m$ rational independent irrational numbers.
Proposition \ref{prop:3} and the previous discussion indicate that
$$
p=\mathcal{N}_{\rho(L_g^d)}=\mathcal{N}_{\rho(L_{g'}^d)}=p'.
$$
Moreover, $\overline{Orb_{L_g^d}(e)}\,\cong\,\overline{Orb_{L_{g'}^d}(e)}$ is homeomorphic the disjoint union of a group of $T^m$, and the number of these mutually disjoint $T^m$ is just $p$.
So we have
$$
h(\overline{Orb_{L_g^d}(e)})=\overline{Orb_{L_{g'}^d}(e)}\cong\mathbb{Z}_p \times T^m,
$$
where $\mathbb{Z}_p$ is a subgroup of ${\rm SU}(2)$, and $T^m$ is naturally a subgroup $T^n$.
It follows from the definition of lens spaces that
$$
{\rm SU}(2) \times T^n/\overline{Orb_{L_g^d}(e)}
\,\cong\,{\rm SU}(2)/\mathbb{Z}_p \times T^n/T^m
\,\cong\,L(p, -1) \times T^{n-m}.
$$
According to Lemma \ref{lem:5}, one can see that there exists some $h':\,\,L(p, -1) \times T^{n-m} \rightarrow \,L(p, -1) \times T^{n-m}$ induced by the topological conjugacy $h$ such that
$$
\pi_0 \circ h=h' \circ \pi_0,
$$
where $\pi_0:\,{\rm SU}(2) \times T^n \rightarrow L(p, -1) \times T^{n-m}$ is the quotient map.
Investigate the following diagram.
$$
\xymatrixcolsep{2pc}
\xymatrix{
{\,\,{{\rm SU}}(2)\,} \ar[d]_-{\pi_0'} \ar[r]^-{i}
& {\,\,{{\rm SU}}(2) \times T^n\,\,} \ar[d]_-{\pi_0} \ar[r]^-{h}
& {\,\,{{\rm SU}}(2) \times T^n\,\,} \ar[d]_-{\pi_0} \ar[r]^-{\pi}
& {\,\,{{\rm SU}}(2)\,\,} \ar[d]_-{\pi_0'}\\
{\,\,L(p, -1)\,\,} \ar[r]_-{i'}
& {\,\,L(p, -1) \times T^{n-m}\,\,} \ar[r]_-{h'}
& {\,\,L(p, -1) \times T^{n-m}\,\,} \ar[r]_-{\pi'}
& {\,\,L(p, -1)\,\,}.}
$$
In the diagram, $\pi_0$ is the quotient map, $\pi_0'$ is the universal covering map, $i,\,i'$ are natural inclusions, $\pi,\,\pi'$ are projections,
and $i,\,i',\,\pi,\,\pi'$ satisfy
$$
\pi_0 \circ i=i' \circ \pi_0',
\qquad \pi_0' \circ \pi=\pi' \circ \pi_0,
$$
Thus, the above diagram is commutative, i.e.,
$$
f' \circ \pi_0'=\pi_0' \circ f,
$$
where $f=\pi \circ h \circ i,\,f'=\pi' \circ h' \circ i'$.
By the definition of the fundamental group of $L(p, -1)$, we obtain a commutative diagram as follows.
$$
\xymatrixcolsep{3pc}
\xymatrix{
{\,\,I\,\,} \ar[r]^{\tilde{\alpha}} \ar[rd]_{\alpha}
& {\,\,{{\rm SU}}(2)\,\,} \ar[r]^{f} \ar[d]^{\pi_0'} &
{\,\,{{\rm SU}}(2)\,\,} \ar[d]^{\pi_0'}\\
&{\,\,L(p, -1)\,\,} \ar[r]^{f'} & {\,\,L(p, -1)\,\,},}
$$
In this diagram, $I=[0, 1]$, $\pi_0',\,f$ and $f'$ are defined above,
$\alpha$ is a loop in $L(p, -1)$ with the base point $\pi'_0(e_0)$,
where $e_0=\left(\begin{array}{cc}
1&0\\
0&1
\end{array}\right)$ is the identify element of ${\rm SU}(2)$,
$\tilde{\alpha}$ is the unique lift of $\alpha$ satisfying $\tilde{\alpha}(0)=e_0$.
In this way, we denote the fundamental group of $L(p, -1)$ by
$$
\pi_1(L(p,-1))=\{[\alpha_j];\,\,j=1,\,2,\,\cdots,\,p\} \cong \mathbb{Z}_p,
$$
where $\alpha_j$ is a loop in $L(p, -1)$ with base point $\pi_0'(e_0)$,
and $[\alpha_j]$ represents the homotopy class of $\alpha_j$.
Let $\tilde{\alpha}_j$ be the unique lift of $\alpha_j$ satisfying $\tilde{\alpha_j}(0)=e_0$.
Then one can see that $\tilde{\alpha}_j$ is a continuous curve  in ${\rm SU}(2)$ connecting the elements $e_0=\left(\begin{array}{cc}
1&0\\
0&1
\end{array}\right)$ and $\left(\begin{array}{cc}
{{\textrm{e}}}^{2\pi{\rm i} j/p}&0\\
0&{{\textrm{e}}}^{-2\pi{\rm i} j/p}
\end{array}\right)$.
According to the conditions
$$
h \circ L_g^d=L_{g'}^d \circ h,\qquad h(e)=e,\qquad f=\pi \circ h \circ i,
$$
it is easy to prove that
$$
f(\tilde{\alpha}_q)=\tilde{\alpha}_{q'},
\qquad \tilde{\alpha}_q,\,\tilde{\alpha}_{q'} \subseteq {\rm SU}(2).
$$
Moreover, these two commutative diagrams indicate that
$$
f'(\alpha_q)=\alpha_{q'},\quad\hbox{i.e.,}\quad f'_*([\alpha_q])=[\alpha_{q'}],
$$
where $f'_*$ is the endomorphism of $\pi_1(L(p, -1))$ induced by $f'$.
On the other hand, it follows from Lemma \ref{lem:3} and \cite{SH-2010} that $f'$ is homotopic to some preserving-orientation self-homeomorphism on $L(p, -1)$.
And then Lemma \ref{lem:4} tells us
$$
f'_*([\alpha_q])=\pm\,[\alpha_q]=[\alpha_{q'}].
$$
So
$$
q'=q\quad({\rm mod}\,\,\,p) \quad\hbox{or}\quad q+q'=0\quad({\rm mod}\,\,\,p),
$$
i.e.,
$$
\dfrac{q'}{p}=\pm\,\dfrac{q}{p}\quad(\rm{mod}\,\,\,\mathbb{Z}),
$$
and then
$$
d\theta'_0=\pm\,d\theta_0+N,\qquad N \in \mathbb{Z}.
$$
Assume that $N = N'\,\,\,(\rm{mod}\,\,\,d)$.
Since ${\rm gcd}\,(d, k)=1$, then there exists some $N'' \in \mathbb{Z}$ such that $N''k = d-N'\,\,\,({\rm mod}\,\,\,d)$.
Thus,
$$
\theta'_0=\pm\,\theta_0+\dfrac{N}{d}-N''\beta_{m+1}+N''\beta_{m+1}=\pm\,\theta_0+N''\beta_{m+1}\quad({\rm mod}\,\,\,\mathbb{Z}),
$$
where $\beta_{m+1}=\dfrac{k}{d}$.
Take
$$
\theta_i=\beta_i,\,\,\,i=1,\,2,\,\cdots,\,m+1,
\qquad\theta_{i'}=0,\,\,\,i'=m+2,\,m+3,\,\cdots,\,n,
$$
and
$$
\ell_j=0,\,\,\,j=1,\,2,\,\cdots,\,m,
\quad \ell_{m+1}=N'',
\quad \ell_{j'} \in \mathbb{Z},\,\,\,j'=m+2,\,m+3,\,\cdots,\,n.
$$
Consequently, we obtain
$$
\theta'_0=\pm\,\theta_0+\sum_{i=1}^n{\ell_i\theta_i}\quad({\rm{mod}}\,\,\,\mathbb{Z}),
\qquad \ell_i \in \mathbb{Z},\,\,\, i=1,\,2,\,\cdots,\,n.
$$

$\bm {(2)}$ Assume that both of $\rho(L_g)$ and $\rho(L_{g'})$ satisfy the conditions in {\bf Case 2}.

Let $h$ be the topological conjugacy from $L_g$ to $L_{g'}$ preserving the identify element.
Then we have
$$
h \circ L_g^d=L_{g'}^d \circ h,\qquad h(\overline{Orb_{L_g^d}(e)})=\overline{Orb_{L_{g'}^d}(e)}.
$$
i.e., $L_g^d$ and $L_{g'}^d$ are topologically conjugate.
Obviously, $L_g^d,\,L_{g'}^d \in \mathscr{M}_{T_{{\rm SU}(2) \times T^n}}$ satisfy
$$
\rho(L_g^d)=\left(\begin{array}{c}
d\theta_0\\
d\beta_1\\
\vdots\\
d\beta_m\\
k\\
0\\
\vdots\\
0
\end{array}\right)\quad({\rm mod}\,\,\,\mathbb{Z}^{n+1}),
\qquad \rho(L_{g'}^d)=\left(\begin{array}{c}
d\theta'_0\\
d\beta'_1\\
\vdots\\
d\beta'_m\\
k'\\
0\\
\vdots\\
0
\end{array}\right)\quad({\rm mod}\,\,\,\mathbb{Z}^{n+1}),
$$
where $d\theta_0,\,d\beta_1,\,d\beta_2,\,\cdots\,d\beta_m$ and $d\theta'_0,\,d\beta'_1,\,d\beta'_2,\,\cdots,\,d\beta'_m$ are both rational independent.
Thus, we get
$$
\widetilde{\mathcal{R}}_{\rho(L_g^d)}=\widetilde{\mathcal{R}}_{\rho(L_{g'}^d)}=m+1
$$
And then Proposition \ref{prop:3} indicates that
$$
\overline{Orb_{L_g^d}(e)}\,\cong\,\overline{Orb_{L_{g'}^d}(e)}\,\cong\,T^{m+1}.
$$
Since $L^d_g,\,L^d_{g'} \in \mathscr{M}_{T_{{\rm SU}(2) \times T^n}}$ and $e \in T_{{\rm SU}(2) \times T^n}$,
then
$$
\overline{Orb_{L_g^d}(e)},\quad\overline{Orb_{L_{g'}^d}(e)} \subseteq T_{{\rm SU}(2) \times T^n}.
$$
In fact, by the forms of $\rho(L_g^d)$ and $\rho(L_{g'}^d)$, one can see that
$$
\overline{Orb_{L_g^d}(e)}=\overline{Orb_{L_{g'}^d}(e)} \subseteq T_{{\rm SU}(2) \times T^n}.
$$
Naturally, we can restrict this dynamical system to $\overline{Orb_{L_g^d}(e)}$, and then obtain
$$
h|_{\overline{Orb_{L_g^d}(e)}} \circ L_g^d|_{\overline{Orb_{L_g^d}(e)}}
=L_{g'}^d|_{\overline{Orb_{L_g^d}(e)}} \circ h|_{\overline{Orb_{L_g^d}(e)}},
$$
where $h|_{\overline{Orb_{L_g^d}(e)}}$ is a self-homeomorphism of $\overline{Orb_{L_g^d}(e)}$,
that means $L_g^d|_{\overline{Orb_{L_g^d}(e)}}$ and $L_{g'}^d|_{\overline{Orb_{L_g^d}(e)}}$ are topologically conjugate.
According to the previous discussion at the beginning of this section,
we know that $L_g^d|_{\overline{Orb_{L_g^d}(e)}}$ and $L_{g'}^d|_{\overline{Orb_{L_g^d}(e)}}$ are topologically conjugate to some rotations $f$ and $f'$ of $T^{m+1}$, respectively, with
$$
\rho(f)=\left(\begin{array}{c}
d\theta_0\\
d\beta_1\\
\vdots\\
d\beta_m\\
\end{array}\right)\quad({\rm mod}\,\,\,\mathbb{Z}^{m+1}),
\qquad \rho(f')=\left(\begin{array}{c}
d\theta'_0\\
d\beta'_1\\
\vdots\\
d\beta'_m
\end{array}\right)\quad({\rm mod}\,\,\,\mathbb{Z}^{m+1}),
$$
and then $f$ and $f'$ are topologically conjugate.
So
$$
\rho(f')={\bm A}'\rho(f)\quad({\rm mod}\,\,\,\mathbb{Z}^{m+1}),
\qquad {\bm A}' \in {\rm GL}_{m+1}(\mathbb{Z}).
$$
Take two rotations $f_0,\,f_0'$ of $T^{n+1}$ with
$$
\rho(f_0)=\rho(L_g^d),
\qquad \rho(f_0')=\rho(L_{g'}^d).
$$
Then we have
$$
\rho(f_0')=\left(\begin{array}{cc}
{\bm A}'&{\bm 0}\\
{\bm 0}&{\bm I}_{n-m}
\end{array}\right)\rho(f_0)\quad({\rm mod}\,\,\,\mathbb{Z}^{n+1}),
$$
where ${\bm I}_{n-m}$ is the identify matrix with order $n-m$.
Thus, $f_0$ and $f_0'$ are topologically conjugate.
Then there exists some ${\bm A}'' \in {\rm GL}_{n+1}(\mathbb{Z})$ such that
$$
\rho(f_0')={\bm A}''\rho(f_0)\quad({\rm mod}\,\,\,\mathbb{Z}^{n+1}),
\quad\hbox{i.e.,}\quad
\rho(L_{g'}^d)={\bm A}''\rho(L_g^d)\quad({\rm mod}\,\,\,\mathbb{Z}^{n+1}).
$$
And the previous analysis illuminates that
$$
\left(\begin{array}{c}
d\beta'_1\\
\vdots\\
d\beta'_m\\
k'\\
0\\
\vdots\\
0
\end{array}\right)=\left(\begin{array}{cc}
{\bm 0}&{\bm A}
\end{array}\right)\left(\begin{array}{c}
d\theta_0\\
d\beta_1\\
\vdots\\
d\beta_m\\
k\\
0\\
\vdots\\
0
\end{array}\right)\quad({\rm mod}\,\,\,\mathbb{Z}^n),
\qquad {\bm A} \in {\rm GL}_{n}(\mathbb{Z}).
$$
It is easy to verify ${\bm A}''=\left(\begin{array}{cc}
\pm\,1&{\bm \hbar}\\
{\bm 0}&{\bm A}
\end{array}\right)$, where ${\bm \hbar}=(\begin{array}{cccc}
\hbar_1&\hbar_2&\cdots&\hbar_n
\end{array})$ is a $1 \times n$ integer matrix.
And then we have
$$
d\theta'_0=\pm\,d\theta_0+\sum_{i=1}^m{d\hbar_i\beta_i}+N,
\qquad N \in \mathbb{Z},
$$
i.e.,
$$
\theta'_0=\pm\,\theta_0+\sum_{i=1}^m{\hbar_i\beta_i}+\dfrac{N}{d},
\qquad N \in \mathbb{Z}.
$$
Assume that $N = N'\,\,\,({\rm mod}\,\,\,d)$.
Since ${\rm gcd}\,(d, k)=1$, then there exists some $\hbar_{m+1} \in \mathbb{Z}$ such that $\hbar_{m+1} k = N'\,\,\,({\rm mod}\,\,\,d)$.
So
$$
\theta'_0=\pm\,\theta_0+\sum_{i=1}^{m}{\hbar_i\beta_i}+\dfrac{N}{d}-\hbar_{m+1}\beta_{m+1}+\hbar_{m+1}\beta_{m+1}
=\pm\,\theta_0+\sum_{i=1}^{m+1}{\hbar_i\beta_i}\,\,\,({\rm mod}\,\,\,\mathbb{Z}),
$$
where $\beta_{m+1}=\dfrac{k}{d}$.
Take
$$
\theta_i=\beta_i,\,\,\,i=1,\,2,\,\cdots,\,m+1,
\qquad\theta_{i'}=0,\,\,\,i'=m+2,\,m+3,\,\cdots,\,n,
$$
and
$$
\ell_j=\hbar_j,\,\,\,j=1,\,2,\,\cdots,\,m+1
\qquad \ell_{j'} \in \mathbb{Z},\,\,\,j'=m+2,\,m+3,\,\cdots,\,n,
$$
Consequently, we obtain
$$
\theta'_0=\pm\,\theta_0+\sum_{i=1}^n{\ell_i\theta_i}\quad({\rm{mod}}\,\,\,\mathbb{Z}),
\qquad \ell_i \in \mathbb{Z},\,\,\, i=1,\,2,\,\cdots,\,n.
$$

$\bm {(3)}$ Assume that both of $\rho(L_g)$ and $\rho(L_{g'})$ satisfy the conditions in {\bf Case 3}.

Let $h$ be the topological conjugacy from $L_g$ to $L_{g'}$ preserving the identify element.
Then we have
$$
h \circ L_g^d=L_{g'}^d \circ h,\qquad h(e)=e,
$$
i.e., $L_g^d$ and $L_{g'}^d$ are topologically conjugate with
$$
\rho(L_g^d)=\left(\begin{array}{c}
\alpha_0\\
\alpha_1\\
\vdots\\
\alpha_m\\
0\\
\vdots\\
0
\end{array}\right),
\qquad \rho(L_{g'}^d)=\left(\begin{array}{c}
\alpha'_0\\
\alpha'_1\\
\vdots\\
\alpha'_m\\
0\\
\vdots\\
0
\end{array}\right),
$$
where
\begin{gather*}
\alpha_0=d\theta_0\,\,\,({\rm mod}\,\,\,\mathbb{Z}),
\qquad \alpha'_0=d\theta'_0\,\,\,({\rm mod}\,\,\,\mathbb{Z}),
\\
\alpha_i=d\beta_i\,\,\,({\rm mod}\,\,\,\mathbb{Z}),
\qquad \alpha'_i=d\beta'_i\,\,\,({\rm mod}\,\,\,\mathbb{Z}),
\qquad i=1,\,2,\,\cdots,\,m.
\end{gather*}
According to the above discussion, we know
$$
\left(\begin{array}{c}
\alpha'_1\\
\vdots\\
\alpha'_m\\
0\\
\vdots\\
0
\end{array}\right)={\bm A}\left(\begin{array}{c}
\alpha_1\\
\vdots\\
\alpha_m\\
0\\
\vdots\\
0
\end{array}\right)\quad({\rm mod}\,\,\,\mathbb{Z}^n),
\qquad {\bm A} \in {\rm GL}_n(\mathbb{Z}).
$$
Take $L' \in \mathscr{M}_{T_{{\rm SU}(2) \times T^n}}$ with
$$
\rho(L')=(\begin{array}{cccccccc}
\alpha'_0&\alpha_1&\alpha_2&\cdots&\alpha_m&0&\cdots&0
\end{array})^{\rm T}.
$$
Then
$$
\rho(L_{g'}^d)=\left(\begin{array}{cc}
1&{\bm 0}\\
{\bm 0}&{\bm A}
\end{array}\right)\rho(L')\quad({\rm mod}\,\,\,\mathbb{Z}^{n+1}).
$$
So the sufficiency of Theorem \ref{the:1} tell us that $L_{g'}^d$ and $L'$ are topologically conjugate.
Set $L=L_g^d$.
And then one can see that $L$ and $L'$  are topologically conjugate.

In fact, $\alpha_1,\,\alpha_2,\,\cdots,\,\alpha_m$ are rational independent,
but $\alpha_0,\,\alpha_1,\,\cdots\,\alpha_m$ are $m+1$ rational dependent irrational numbers,
and $\alpha'_0,\,\alpha_1,\,\cdots\,\alpha_m$ are $m+1$ rational dependent irrational numbers, too.
Then we have
\begin{align}\label{eq:5}
\left\{
\begin{array}{l}
k_0\alpha_0+k_1\alpha_1+\cdots+k_m\alpha_m+\dfrac{q}{p}=0,
\\[0.3 cm]
k'_0\alpha'_0+k'_1\alpha_1+\cdots+k'_m\alpha_m+\dfrac{q'}{p'}=0,
\end{array}
\right.
\end{align}
where $k_0,\,k'_0$ are positive integers, $p,\,p',\,q,\,q',\,k_i,\,k_i' \in \mathbb{Z},\,i=1,\,2,\cdots,\,m$, and
$$
{\rm gcd}\,(k_0, k_1, \cdots, k_m)={\rm gcd}\,(k'_0, k'_1, \cdots, k'_m)={\rm gcd}\,(p, q)={\rm gcd}\,(p', q')=1.
$$
It follows from Proposition \ref{prop:3} and the previous discussion about topologically conjugate invariants that we have
$$
\widetilde{\mathcal{R}}_{\rho(L)}=\widetilde{\mathcal{R}}_{\rho(L')}=m,
\qquad \mathcal{N}_{\rho(L)}=p,
\qquad \mathcal{N}_{\rho(L')}=p'.
$$
And then $L$ and $L'$ being topologically conjugate implies that $p'=p$.
Let $h'$ be the topologically conjugate form $L$ to $L'$ satisfying $h'(e)=e$.
Then
$$
L^p,\,L'^p \in \mathscr{M}_{T_{{\rm SU}(2) \times T^n}},
\qquad h' \circ L^p=L'^p \circ h',
$$
and
$$
\rho(L^p)=\left(\begin{array}{c}
p\alpha_0\\
p\alpha_1\\
\vdots\\
p\alpha_m\\
0\\
0\\
\vdots\\
0
\end{array}\right)\quad({\rm mod}\,\,\,\mathbb{Z}^{n+1}),
\qquad \rho(L'^p)=\left(\begin{array}{c}
p\alpha'_0\\
p\alpha_1\\
\vdots\\
p\alpha_m\\
0\\
0\\
\vdots\\
0
\end{array}\right)\quad({\rm mod}\,\,\,\mathbb{Z}^{n+1}).
$$
Lemma \ref{lem:6} indicates that there exist some matrices ${\bm G},\,{\bm G}' \in {\rm GL}_{m+1}(\mathbb{Z})$ such that $(k_0, k_1, \cdots, k_m)$ and $(k'_0, k'_1, \cdots, k'_m)$ are just the first rows of ${\bm G}$ and ${\bm G}'$, respectively.
Let $F,\,F':\,T^m \rightarrow S
^1 \times T^m$ be two embedding maps satisfying
$$
F(e')=e'',\qquad F'(e')=e'',
$$
and
$$
F_*({\bm l})={\bm G}_1{\bm l},\qquad F'_*({\bm l})={\bm G}'_1{\bm l},
\qquad \forall\,{\bm l} \in \pi_1(T^m)\,\cong\,\mathbb{Z}^m,
$$
where $e',\,e''$ are the identity elements of $T^m$ and $S^1 \times T^m\,\cong\,T^{m+1}$, respectively,
$F_*,\,F'_*:\,\pi_1(T^m) \rightarrow \pi_1(T^{m+1})$ are the group homomorphisms induced by $F$ and $F'$, respectively,
and ${\bm G}_1,\,{\bm G}'_1$ are the $(m+1) \times m$ integer matrices obtained by removing the first columns of ${\bm G}^{-1}$ and ${\bm G}'^{-1}$, respectively, where ${\bm G}^{-1},\,{\bm G}'^{-1} \in {\rm GL}_{m+1}(\mathbb{Z})$.
Investigate the following diagram.
$$
\xymatrixcolsep{3pc}
\xymatrix{
{\,\,T^m\,} \ar[d]_-{h''} \ar[r]^-{F}
& {\,\,S^1 \times T^m\,\,}  \ar[r]^-{i}
& {\,\,{\rm SU}(2) \times T^n\,\,} \ar[d]_-{h'} \\
{\,\,T^m\,\,} \ar[r]_-{F'}
& {\,\,S^1 \times T^m\,\,} \ar[r]_-{i}
& {\,\,{\rm SU}(2) \times T^n\,\,}.}
$$
In the diagram, $F,\,F'$ are the embedding maps defined above,
$h'$ is the topological conjugacy from $L^p$ to $L'^p$ satisfying $h(e)=e$,
$h''$ is a self-homeomorphism of $T^m$,
$i$ is the natural inclusion preserving the identity elements and satisfying
$$
i(S^1) \subseteq {\rm SU}(2),\qquad i(T^m) \subseteq T^n.
$$
Thus, it is easy to know that
$$
i_*({\bm l}')=\left(\begin{array}{cc}
{\bm 0}&{\bm I}_m\\
{\bm 0}&{\bm 0}
\end{array}\right){\bm l}'=\left(\begin{array}{c}
l_1\\
l_2\\
\vdots\\
l_m\\
0\\
\vdots\\
0
\end{array}\right),
\qquad \forall\,{\bm l}'=\left(\begin{array}{c}
l_0\\
l_1\\
\vdots\\
l_m\\
\end{array}\right) \in \pi_1(T^{m+1})\,\cong\,\mathbb{Z}^{m+1},
$$
where $i_*$ is the group homomorphism between fundamental groups induce by $i$,
$\left(\begin{array}{cc}
{\bm 0}&{\bm I}_m\\
{\bm 0}&{\bm 0}
\end{array}\right)$ is an $n \times (m+1)$ integer matrix, ${\bm I}_m$ is the identity matrix with order $m$.
Set
$$
H=i \circ F:\,T^m \rightarrow i \circ F(T^m)\,\cong\,T^m,\quad H'=i \circ F':\,T^m \rightarrow i \circ F'(T^m)\,\cong\,T^m.
$$
Then $H$ $H'$ are both homeomorphisms satisfying
$$
H(e')=e,\qquad H'(e')=e,
$$
where $e,\,e'$ are the identity elements of ${\rm SU}(2) \times T^n$ and $T^m$, respectively.
And take two rotations of $f,\,f'$ of $T^m$ with
$$
\rho(f)={\bm G}_2\left(\begin{array}{c}
p\alpha_0\\
p\alpha_1\\
\vdots\\
p\alpha_m
\end{array}\right)\quad({\rm mod}\,\,\,\mathbb{Z}^m),
\quad \rho(f')={\bm G}'_2\left(\begin{array}{c}
p\alpha'_0\\
p\alpha_1\\
\vdots\\
p\alpha_m
\end{array}\right)\quad({\rm mod}\,\,\,\mathbb{Z}^m),
$$
where ${\bm G}_2,\,{\bm G}'_2$ are the $m \times (m+1)$ integer matrices obtained by removing the first rows of ${\bm G}$ and ${\bm G}'$, respectively.
It follows from Proposition \ref{prop:1} that
$$
\widetilde{\mathcal{R}}_{\rho(f)}=\widetilde{\mathcal{R}}_{\rho(L)}=m,
\qquad \widetilde{\mathcal{R}}_{\rho(f')}=\widetilde{\mathcal{R}}_{\rho(L')}=m.
$$
Then $\overline{Orb_f(e')}=\overline{Orb_{f'}(e')}=T^m$.
Thus, through a simple calculation, we have
$$
H \circ f=L^p|_{H(T^m)} \circ H,
\qquad H' \circ f'=L'^p|_{H'(T^m)} \circ H',
$$
and then
\begin{gather*}
H(\overline{Orb_f(e')})=\overline{Orb_{L^p|_{H(T^m)}}(e)}=\overline{Orb_{L^p}(e)}\,\cong\,T^m,
\\
H'(\overline{Orb_{f'}(e')})=\overline{Orb_{L'^p|_{H'(T^m)}}(e)}=\overline{Orb_{L'^p}(e)}\,\cong\,T^m.
\end{gather*}
That means
$$
H(T^m)=\overline{Orb_{L^p}(e)},
\qquad H'(T^m)=\overline{Orb_{L'^p}(e)}.
$$
Similar to the discussion in Case $\bm{(2)}$,
we can restrict the dynamical system to $\overline{Orb_{L^p}(e)} \cong \overline{Orb_{L'^p}(e)}$,
and then obtain
$$
h'|_{\overline{Orb_{L^p}(e)}} \circ L^p|_{\overline{Orb_{L^p}(e)}}=L'^p|_{\overline{Orb_{L'^p}(e)}} \circ h'|_{\overline{Orb_{L^p}(e)}},
$$
i.e,
$$
h'|_{H(T^m)} \circ L^p|_{H(T^m)}=L'^p|_{H'(T^m)} \circ h'|_{H(T^m)}.
$$
Consequently, $f$ and $f'$ are topologically conjugate.
Set $h''=H'^{-1} \circ h' \circ H$.
And it is not difficult to verify that
$$
h'' \circ f=f' \circ h'',
$$
i.e., $h''$ is a topological conjugate from $f$ to $f'$.
And then according to the discussion in introduction, we have
$$
\rho(f')={\bm A}'\rho(f)\quad({\rm mod}\,\,\,\mathbb{Z}),
\qquad {\bm A}' \in {\rm GL}_m(\mathbb{Z}),
$$
where ${\bm A}' \in {\rm GL}_m(\mathbb{Z})$ is just the matrix form of the automorphism $h''_*$ of $\pi_1(T^m)$ induced by $h''$.
On the other hand, the condition $h''=H'^{-1} \circ h' \circ H$ indicates that the above diagram is commutative.
So we have
$$
i\circ F' \circ h''=h' \circ i \circ F,
$$
and then
$$
i_* \circ F'_* \circ h''_*=h'_* \circ i_* \circ F_*,
$$
where $h'_*$ is the automorphism of $\pi_1({\rm SU}(2) \times T^n)$ induced by $h'$ satisfying
$$
h'_*({\bm l}'')={\bm A}''{\bm l}'',\qquad \forall\,{\bm l}'' \in \pi_1({\rm SU}(2) \times T^n)\,\cong\,\mathbb{Z}^n,
$$
where ${\bm A}'' \in {\rm GL}_n(\mathbb{Z})$ is just the matrix form of $h'_*$.
Thus, for any ${\bm l} \in \pi_1(T^m)\,\cong\,\mathbb{Z}^m$, we get
$$
\left(\begin{array}{cc}
{\bm 0}&{\bm I}_m\\
{\bm 0}&{\bm 0}
\end{array}\right){\bm G}'_1{\bm A}'({\bm l})={\bm A}''\left(\begin{array}{cc}
{\bm 0}&{\bm I}_m\\
{\bm 0}&{\bm 0}
\end{array}\right){\bm G}_1({\bm l}),
$$
i.e.,
$$
\left(\begin{array}{c}
{\bm K'}\\
{\bm 0}
\end{array}\right){\bm A}'={\bm A}''\left(\begin{array}{c}
{\bm K}\\
{\bm 0}
\end{array}\right),
$$
where $\left(\begin{array}{c}
{\bm K}\\
{\bm 0}
\end{array}\right)$ and $\left(\begin{array}{c}
{\bm K'}\\
{\bm 0}
\end{array}\right)$ are $n \times m$ integer matrices,
${\bm K},\,{\bm K'}$ are the integer matrices with order $m$ obtained by removing the first rows of ${\bm G}_1$ and ${\bm G}'_1$, respectively.
Obviously, we can directly obtain ${\bm K},\,{\bm K'}$ by removing the first rows and first columns of ${\bm G}^{-1}$ and ${\bm G}'^{-1}$, respectively.
Moreover, we know that the element at the first row and first column of ${\bm G}$ is $k_0$, and for ${\bm G}'$, the element at the same position is $k_0'$.
So it follows from the properties of invertible matrices that
$$
|{\bm K}|=k_0,\qquad |{\bm K'}|=k'_0.
$$
In fact, considering the forms of $\rho(L^p)$ and $\rho(L'^p)$, we can write ${\bm A}''$ by
$$
{\bm A}''=\left(\begin{array}{cc}
{\bm I}_m&{\bm B}'\\
{\bm 0}&{\bm B}
\end{array}\right),
$$
where ${\bm B} \in {\rm GL}_{n-m}(\mathbb{Z})$,
${\bm B}'$ is an $m \times (n-m)$ integer matrix.
Thus, through a simple verification, one can see that ${\bm K'}{\bm A}'={\bm K}$,
that means $|{\bm K'}||{\bm A}'|=|{\bm K}|$.
And then combining this fact with the known conditions $|{\bm A}'|=\pm\,1$ and $k_0,\,k'_0 \in \mathbb{Z}_+$,
we have
$$
k'_0=|{\bm K'}|=|{\bm K}|=k_0.
$$
Consequently, we can write \eqref{eq:5} by
\begin{align}\label{eq:6}
\left\{
\begin{array}{l}
\alpha_0+k_1\dfrac{\alpha_1}{k_0}+\cdots+k_m\dfrac{\alpha_m}{k_0}+\dfrac{q}{pk_0}=0,
\\[0.3 cm]
\alpha'_0+k'_1\dfrac{\alpha_1}{k_0}+\cdots+k'_m\dfrac{\alpha_m}{k_0}+\dfrac{q'}{pk_0}=0.
\end{array}
\right.
\end{align}

Regard ${\rm SU}(2) \times T^n$ as a $k_0^m$-fold covering space of itself,
and define the covering map $\pi:\,{\rm SU}(2) \times T^n \rightarrow {\rm SU}(2) \times T^n$ by
$$
\pi:\,(u, z_1, z_2, \cdots, z_n) \mapsto (u, z_1^{k_0}, z_2^{k_0}, \cdots, z_m^{k_0}, z_{m+1}, \cdots, z_n),
$$
where $u \in {\rm SU}(2)$, $z_i \in \mathbb{C}$ and $|z_i|=1,\,i=1,\,2,\,\cdots,\,n$.
According to the map lifting theorem, one can see that there always exist some lifts of the self-homeomorphisms of ${\rm SU}(2) \times T^n$ under the the covering map $\pi$.
Fix one lift $\tilde{L}_0$ of $L$ and one lift $\tilde{h'}$ of $h'$ satisfying
$$
\rho(\tilde{L}_0)=\left(\begin{array}{ccccccc}
\alpha_0&\alpha_1/k_0&\cdots&\alpha_m/k_0&0&\cdots&0
\end{array}\right)^{\rm T},\qquad \tilde{h'}(e)=e.
$$
Then Lemma \ref{lem:7} indicates that there exists some lift $\tilde{L}'_0$ of $L'$ such that
$$
\tilde{h'} \circ \tilde{L}_0=\tilde{L}'_0 \circ \tilde{h'},
$$
i.e., $\tilde{h'}$ is a topological conjugacy from $\tilde{L}_0$ to $\tilde{L}'_0$.
And it is not difficult to prove that $h'_*=\tilde{h'}_*$ by the properties of covering spaces and covering maps,
where $h'_*,\,\tilde{h'}_*$ are the automorphisms of $\pi_1({\rm SU}(2) \times T^n)$ induced by $h'$ and $\tilde{h'}$, respectively.
If we set
$$
\rho(\tilde{L}'_0)=\left(\begin{array}{c}
\alpha'_0\\
\alpha_1/k_0+p_1/k_0\\
\vdots\\
\alpha_m/k_0+p_m/k_0\\
0\\
\vdots\\
0
\end{array}\right)\,({\rm mod}\,\mathbb{Z}^{n+1}),
\,\,\, p_1,\,\cdots,\,p_m \in \{0, 1, \cdots, k_0-1\},
$$
then
$$
\left(\begin{array}{c}
\alpha_1/k_0+p_1/k_0\\
\alpha_2/k_0+p_2/k_0\\
\vdots\\
\alpha_m/k_0+p_m/k_0\\
0\\
\vdots\\
0
\end{array}\right)={\bm A}''\left(\begin{array}{c}
\alpha_1/k_0\\
\alpha_2/k_0\\
\vdots\\
\alpha_m/k_0\\
0\\
\vdots\\
0
\end{array}\right)\quad({\rm mod}\,\,\,\mathbb{Z}^n),
$$
where
$$
{\bm A}''=\left(\begin{array}{cc}
{\bm I}_m&{\bm B}'\\
{\bm 0}&{\bm B}
\end{array}\right)
$$
is the matrix form of $h'_*=\tilde{h'}_*$ defined in the previous discussion.
Thus, it is obvious that $p_i=0,\,i=1,\,2,\,\cdots,\,m$, that means
$$
\rho(\tilde{L}'_0)=\left(\begin{array}{ccccccc}
\alpha'_0&\alpha_1/k_0&\cdots&\alpha_m/k_0&0&\cdots&0
\end{array}\right)^{\rm T}.
$$
Set
$$
G=\left\{\left(\begin{array}{cc}
z_0z_1 \cdots z_n&0\\
0&\bar{z_0}
\end{array}\right);\,z_0z_1^{k_1} \cdots z_m^{k_m}=1,\, z_i=1,\,i=m+1,\,\cdots,\,n \right\},
$$
$$
G'=\left\{\left(\begin{array}{cc}
\omega_0\omega_1 \cdots \omega_n&0\\
0&\bar{\omega_0}
\end{array}\right);\,\omega_0\omega_1^{k_1'} \cdots \omega_m^{k_m'}=1,\, \omega_i=1,\,i=m+1,\,\cdots,\,n \right\},
$$
where $z_j,\,\omega_j \in \mathbb{C},\,|z_j|=|\omega_j|=1,\,j=0,\,1,\,\cdots,\,m$.
Then according to Lemma \ref{lem:8},
one can see that $G,\,G'$ are two subgroups of ${\rm SU}(2) \times T^n$ satisfying $G\,\cong\,G'\,\cong\,T^m$,
and
$$
{\rm SU}(2) \times T^n/G\,\cong\,{\rm SU}(2) \times T^n/G'
\,\cong\,L(1, -1) \times T^{n-m}
\,\cong\,{\rm SU}(2)\times T^{n-m}.
$$
And then using the above observations, we have
$$
\dfrac{pk_0}{{\rm gcd}\,(k_0, q)}=\mathcal{N}_{\rho(\tilde{L}_0)}=\mathcal{N}_{\rho(\tilde{L}'_0)}=\dfrac{pk_0}{{\rm gcd}\,(k_0, q')},
$$
i.e., ${\rm gcd}\,(k_0, q)={\rm gcd}\,(k_0, q')$. Set
$$
p_0=\dfrac{pk_0}{{\rm gcd}\,(k_0, q)}=\dfrac{pk_0}{{\rm gcd}\,(k_0, q')}.
$$
In fact, it follows from a simple verification that
$$
\overline{Orb_{\tilde{L}^{p_0}_0}(e)}=G,
\qquad\overline{Orb_{\tilde{L}'^{p_0}_0}(e)}=G',
$$
and then $\tilde{h'}(G)=G'$.
So Lemma \ref{lem:5} tell us that there exists some homeomorphism $h_0:\,{\rm SU}(2) \times T^{n-m} \rightarrow {\rm SU}(2) \times T^{n-m}$ such that
$$
\pi_0' \circ \tilde{h'}=h_0 \circ \pi_0,
$$
where
$$
\pi_0:\,{\rm SU}(2) \times T^m \rightarrow {\rm SU}(2) \times T^m/G\,\cong\,{\rm SU}(2) \times T^{n-m},
$$
$$
\pi'_0:\,{\rm SU}(2) \times T^m \rightarrow {\rm SU}(2) \times T^m/G'\,\cong\,{\rm SU}(2) \times T^{n-m}
$$
are quotient maps.
By the same way in the proof of Lemma \ref{lem:8},
we denote the element $u \in {\rm SU}(2) \times T^n$ by
$$
u=(z, {\textrm{e}}^{2\pi{\rm i}\gamma_0}, {\textrm{e}}^{2\pi{\rm i}\gamma_1}, \cdots, {\textrm{e}}^{2\pi{\rm i}\gamma_n}),
$$
where $(z, {\textrm{e}}^{2\pi{\rm i}\gamma_0}) \in {\rm SU}(2),\,\gamma_i \in \mathbb{C},\,i=0,\,1,\,\cdots,\,n$,
and prove that for any
$$
u=(z, {\textrm{e}}^{2\pi{\rm i}\gamma_0}, {\textrm{e}}^{2\pi{\rm i}\gamma_1}, \cdots, {\textrm{e}}^{2\pi{\rm i}\gamma_n}) \in {\rm SU}(2) \times T^n,
$$
the quotient maps $\pi_0$ and $\pi'_0$ satisfy
\begin{gather*}
\pi_0:\,u \mapsto (z, {\textrm{e}}^{2\pi{\rm i}(\gamma_0+k_1\gamma_1+k_2\gamma_2+\cdots+k_m\gamma_m)}, {\textrm{e}}^{2\pi{\rm i}\gamma_{m+1}}, \cdots, {\textrm{e}}^{2\pi{\rm i}\gamma_n}) \in {\rm SU}(2) \times T^{n-m},
\\
\pi'_0:\,u \mapsto (z, {\textrm{e}}^{2\pi{\rm i}(\gamma_0+k'_1\gamma_1+k'_2\gamma_2+\cdots+k'_m\gamma_m)}, {\textrm{e}}^{2\pi{\rm i}\gamma_{m+1}}, \cdots, {\textrm{e}}^{2\pi{\rm i}\gamma_n}) \in {\rm SU}(2) \times T^{n-m},
\end{gather*}
respectively.
Take two left actions $L_0,\,L'_0 \in {\mathscr M}_{T_{{\rm SU}(2) \times T^{n-m}}}$ with
$$
\rho(L_0)=\left(\begin{array}{c}
-q/pk_0\\
0\\
\vdots\\
0
\end{array}\right)\,({\rm mod}\,\mathbb{Z}^{n-m+1}),
\,\rho(L'_0)=\left(\begin{array}{c}
-q'/pk_0\\
0\\
\vdots\\
0
\end{array}\right)\,({\rm mod}\,\mathbb{Z}^{n-m+1}).
$$
Through a simple calculation, one can see that
$$
\pi_0 \circ {\tilde L}_0(u)=L_0 \circ \pi_0(u),\qquad \pi'_0 \circ {\tilde L}'_0(u)=L'_0 \circ \pi'_0(u),
\qquad \forall\,u \in {\rm SU}(2) \times T^{n-m},
$$
Thus, for any $u \in {\rm SU}(2) \times T^n$, we have
\begin{gather*}
\pi'_0 \circ \tilde{h'} \circ {\tilde L}_0(u)=
h_0 \circ \pi_0 \circ {\tilde L}_0(u)=h_0 \circ L_0 \circ \pi_0(u),
\\
\pi'_0 \circ {\tilde L}'_0 \circ \tilde{h'}(u)=
L'_0 \circ \pi'_0 \circ \tilde{h'}(u)=L'_0 \circ h_0 \circ \pi_0(u).
\end{gather*}
And then the known condition $\tilde{h'} \circ {\tilde L}_0={\tilde L}'_0 \circ \tilde{h'}$ implies that
$$
h_0 \circ L_0 \circ \pi_0(u)=L'_0 \circ h_0 \circ \pi_0(u),
\qquad \forall\,u \in {\rm SU}(2) \times T^n.
$$
Obviously, the quotient map $\pi_0:\,{\rm SU}(2) \times T^n \rightarrow {\rm SU}(2) \times T^{n-m}$ is a surjection.
So
$$
h_0 \circ L_0(v)=L'_0 \circ h_0(v),
\qquad \forall\,v \in {\rm SU}(2) \times T^{n-m},
$$
that means $L_0$ and $L_0'$ are topologically conjugate.
Finally, by the conclusions in Case ${\bm {(1)}}$, we have
$$
\dfrac{q'}{pk_0}=\pm\,\dfrac{q}{pk_0}\quad({\rm{mod}}\,\,\mathbb{Z}).
$$
Take another pair of lifts of $L$ and $L'$ denoted by $\tilde{L}_1,\,\tilde{L}'_1$ with
$$
\rho(\tilde{L}_1)=\left(\begin{array}{c}
\alpha_0\\
\alpha_1/k_0+1/k_0\\
\alpha_1/k_0\\
\vdots\\
\alpha_m/k_0\\
0\\
\vdots\\
0
\end{array}\right),
\,\,\, \rho(\tilde{L}'_1)=\left(\begin{array}{c}
\alpha'_0\\
\alpha_1/k_0+1/k_0\\
\alpha_1/k_0\\
\vdots\\
\alpha_m/k_0\\
0\\
\vdots\\
0
\end{array}\right)\,\,\,({\rm mod}\,\,\mathbb{Z}^{n+1}).
$$
It is easy to know that $\tilde{L}_1$ and $\tilde{L}'_1$ are topologically conjugate,
and $\tilde{h'}$ is also a topological conjugacy form $\tilde{L}_1$ to $\tilde{L}'_1$.
Similarly, we can take $L_1,\,L'_1 \in {\mathscr M}_{T_{{\rm SU}(2) \times T^{n-m}}}$ satisfying
$$
\rho(L_1)=\left(\begin{array}{c}
(-q+pk_1)/pk_0\\
0\\
\vdots\\
0
\end{array}\right),
\,\rho(L'_1)=\left(\begin{array}{c}
(-q'+pk'_1)/pk_0\\
0\\
\vdots\\
0
\end{array}\right)\,({\rm mod}\,\mathbb{Z}^{n-m+1}),
$$
and then obtain
$$
h_0 \circ L_1=L_1' \circ h_0.
$$
Consequently, we have
$$
\dfrac{q'-pk'_1}{pk_0}=\pm\,\dfrac{q-pk_1}{pk_0}\quad({\rm{mod}}\,\,\,\mathbb{Z}).
$$
Assume that
$$
q = n_1\quad({\rm{mod}}\,\,\,pk_0),\qquad q' = n'_1\quad({\rm{mod}}\,\,\,pk_0),
$$
and
$$
q-pk_1 = n_2\quad({\rm{mod}}\,\,\,pk_0),\qquad q'-pk_1' = n'_2\quad({\rm{mod}}\,\,\,pk_0),
$$
where $n_1,\,n'_1,\,n_2,\,n'_2 \in \{0,\,1,\,\cdots,\,pk_0-1\}$.

On the other hand, it is known that ${\rm gcd}\,(k_0, k_1, \cdots, k_m)=1$.
Assume that
$$
{\rm gcd}\,(k_1, k_2, \cdots, k_m)=\bar{k}.
$$
So we have ${\rm gcd}\,(k_0, \bar{k})=1$, and
$$
{\rm gcd}\,(k_1/\bar{k},\,\,\,k_2/\bar{k},\,\,\,\cdots,\,\,\,k_m/\bar{k})=1.
$$
Then Lemma \ref{lem:6} illuminates that there exists $l_1,\,l_2,\,\cdots,\,l_m \in \mathbb{Z}$ such that
$$
l_1k_1/\bar{k}+l_2k_2/\bar{k}+\cdots+l_mk_m/\bar{k}=1,
$$
i.e.,
$$
l_1k_1+l_2k_2+\cdots+l_mk_m=\bar{k}.
$$
Thus, take the lifts $\tilde{\mathcal{L}}_i,\,\tilde{\mathcal{L}}'_i$ of $L$ and $L'$, respectively, satisfying
$$
\rho(\tilde{\mathcal{L}}_i)=\left(\begin{array}{c}
\alpha_0\\
\alpha_1/k_0+il_1/k_0\\
\alpha_2/k_0+il_2/k_0\\
\vdots\\
\alpha_m/k_0+il_m/k_0\\
0\\
\vdots\\
0
\end{array}\right),
\,\,\, \rho(\tilde{\mathcal{L}}'_i)=\left(\begin{array}{c}
\alpha'_0\\
\alpha_1/k_0+il_1/k_0\\
\alpha_2/k_0+il_2/k_0\\
\vdots\\
\alpha_m/k_0+il_m/k_0
\\
0\\
\vdots\\
0
\end{array}\right)\,\,\,({\rm mod}\,\,\mathbb{Z}^{n+1}),
$$
and the left actions $\mathcal{L}_i,\,\mathcal{L}'_i \in {\mathscr M}_{T_{{\rm SU}(2) \times T^{n-m}}}$ with
$$
\rho(\mathcal{L}_i)=\left(\begin{array}{c}
-q/pk_0+i\bar{k}/k_0\\
0\\
\vdots\\
0
\end{array}\right),
\,\rho(\mathcal{L}'_i)=\left(\begin{array}{c}
-q'/pk_0+ik'/k_0\\
0\\
\vdots\\
0
\end{array}\right)\,({\rm mod}\,\mathbb{Z}^{n-m+1}),
$$
where $k'=l_1k'_1+l_2k'_2+\cdots+l_mk'_m$.
By the same way, we obtain
$$
\pi \circ \tilde {\mathcal{L}}_i=\mathcal{L}_i \circ \pi,\qquad \pi' \circ \tilde {\mathcal{L}}'_i=\mathcal{L}'_i \circ \pi',
$$
$\tilde{h'}$ is a topological conjugacy from $\tilde{\mathcal{L}}_i$ to $\tilde{\mathcal{L}}'_i$,
and $h_0$ is a topological conjugacy from $\mathcal{L}_i$ to $\mathcal{L}'_i$.
Set $k_0={\rm gcd}\,(k_0, p)\cdot a$. Then $pk_0={\rm gcd}\,(k_0, p)\cdot a \cdot p$.
And the known condition ${\rm gcd}\,(p, q)=1$ implies that ${\rm gcd}\,(-q+ip, p)=1$.
So we have
$$
{\rm gcd}\,(-q+ip,\,\,\,{\rm gcd}\,(k_0, p))=1.
$$
Since ${\rm gcd}\,(a, p)=1$, then there exist some $s,\,t \in \mathbb{Z}$ such that
$$
sp-ta=q+1,\quad\hbox{i.e.}\quad -q+sp=ta+1.
$$
And it is easy to know that ${\rm gcd}\,(-q+sp, a)=1$.
Thus, we get
$$
{\rm gcd}\,(-q+sp,\,\,\,pk_0)=1.
$$
Notice that ${\rm gcd}\,(k_0, \bar{k})=1$, then there exists some $i_0 \in \{0, 1,\cdots, k_0-1\}$ such that $i_0\bar{k} = s\,\,\,({\rm mod}\,\,\,k_0)$,
i.e., ${\rm gcd}\,(-q+i_0\bar{k}p, pk_0)=1$.
We have known that $h_0$ is the topological conjugacy from $\mathcal{L}_{i_0}$ to $\mathcal{L}'_{i_0}$ satisfying $h_0(e_0)=e_0$, where $e_0$ is the identity element of ${\rm SU}(2) \times T^{n-m}$.
Then the fact $(-q+i_0\bar{k}p, pk_0)=1$ indicates that
$$
h_0(\overline{Orb_{\mathcal{L}_{i_0}}(e_0)})=\overline{Orb_{\mathcal{L}'_{i_0}}(e_0)}\,\cong\,\mathbb{Z}_{pk_0},
\,\,\hbox{i.e.,}\,\, \overline{Orb_{\mathcal{L}_{i_0}}(e_0)}\,\cong\,\overline{Orb_{\mathcal{L}'_{i_0}}(e_0)}\,\cong\,\mathbb{Z}_{pk_0},
$$
and then
$$
{\rm SU}(2) \times T^{n-m}/\overline{Orb_{\mathcal{L}_{i_0}}(e_0)}
\,\cong\,{\rm SU}(2)/\mathbb{Z}_{pk_0} \times T^{n-m}\,\cong\,L(pk_0, -1)\times T^{n-m}.
$$
Thus, according to Lemma $\ref{lem:5}$, we know that $h_0$ naturally induced a homeomorphism $h'_0:\,L(pk_0, -1) \rightarrow L(pk_0, -1)$ such that
$$
\pi' \circ h_0=h'_0 \circ \pi',
$$
where $\pi':\,{{\rm SU}}(2) \times T^{n-m} \rightarrow L(pk_0, -1) \times T^{n-m}$ is a quotient map.
Next, we use the same symbols and analysis in Case ${\bm {(1)}}$, and then obtain
$$
f'(\alpha_{n_1})=\alpha_{n_1'},\quad\hbox{i.e.,}\quad f'_*([\alpha_{n_1}])=[\alpha_{n_1'}],
$$
where $f':\,L(pd, -1) \rightarrow L(pd, -1)$ is defined by
$$
\xymatrixcolsep{1.5pc}
\xymatrix{
{L(pd, -1)}  \ar[r]^-{i}
& {L(pd, -1) \times T^{n-m}}  \ar[r]^-{h_0'}
& {L(pd, -1) \times T^{n-m}}  \ar[r]^-{\pi_p}
& {L(pd, -1)},}
$$
$f'_*$ is the endomorphism of $\pi_1(L(pk_0, -1))$ induced by $f'$,
$i$ is the natural inclusion map, $\pi_p$ is the projection.
So Lemma \ref{lem:3} and \cite{SH-2010} indicate that $f'$ is homotopic some orientation-preserving self-homeomorphism on $L(pk_0, -1)$.
And then by Lemma \ref{lem:4}, we have
$$
f'_*([\alpha_{n_1}])=\pm\,[\alpha_{n_1}]=[\alpha_{n'_1}].
$$
Thus, there are two different cases as follows.

$\bm{(a)}$ Assume that $n'_1=n_1$,
i.e., $f'_*([\alpha_{n_1}])=[\alpha_{n_1}]$.
Then Lemma \ref{lem:4} implies that
$$
f'_*(l)=l,
\qquad \forall\,l \in \pi_1(L(pk_0, -1)).
$$
So we get
$$
q' = q\quad({\rm{mod}}\,\,\,pk_0),\quad\hbox{i.e.,}\quad pk_0 \mid q-q'.
$$
And for $L_1,\,L_1'$, one can see that
$$
f'_*([\alpha_{n_2}])=[\alpha_{n_2'}]=[\alpha_{n_2}].
$$
Thus, $n_2'=n_2$, that means
$$
\dfrac{q'-pk'_1}{pk_0}=\dfrac{q-pk_1}{pk_0}\quad({\rm{mod}}\,\,\mathbb{Z}).
$$
Then
$$
k'_1 = k_1\quad({\rm{mod}}\,\,\,k_0),\quad\hbox{i.e.,}\quad k_0 \mid k_1-k'_1.
$$
Moreover, choose the other lifts of $L$ and $L'$, and then by the same way, we can also obtain
$$
k_0 \mid k_i-k'_i,\qquad i=2,\,3,\,\cdots,\,m.
$$
Consequently, combining the above discussion with \eqref{eq:6}, we have
$$
\alpha'_0=\alpha_0+\sum_{i=1}^m\dfrac{k_i-k'_i}{k_0}\alpha_i+\dfrac{q-q'}{pk_0},
$$
i.e.,
$$
\alpha'_0=\alpha_0+\sum_{i=1}^ms_i\alpha_i+N_1,
\qquad s_i,\,N_1 \in \mathbb{Z},\,\,\,i=1,\,2,\,\cdots,\,m.
$$

$\bm{(b)}$ Assume that $n'_1=pk_0-n_1$,
i.e., $f'_*([\alpha_{n_1}])=-[\alpha_{n_1}]$.
Then Lemma \ref{lem:4} implies that
$$
f'_*(l)=-l,
\qquad \forall\,l \in \pi_1(L(pk_0, -1)).
$$
So we get
$$
q+q' = 0\quad({\rm{mod}}\,\,\,pk_0),\quad\hbox{i.e.},\quad pk_0 \mid q+q'.
$$
And for $L_1,\,L_1'$, one can see that
$$
f'_*([\alpha_{n_2}])=[\alpha_{n_2'}]=-[\alpha_{n_2}].
$$
Thus, $n_2'=pk_0-n_2$, that means
$$
\dfrac{q'-pk'_1}{pk_0}=-\dfrac{q-pk_1}{pk_0}\quad({\rm{mod}}\,\,\mathbb{Z}).
$$
Then
$$
k_1+k_1 = 0\quad({\rm{mod}}\,\,\,p),\quad\hbox{i.e.,}\quad k_0 \mid k_1+k'_1.
$$
Similarly, we can also get
$$
k_0 \mid k_i+k'_i,\qquad i=2,\,3,\,\cdots,\,m.
$$
Therefore, \eqref{eq:6} implies that
$$
\alpha'_0=-\alpha_0-\sum_{i=1}^m\dfrac{k_i+k'_i}{k_0}\alpha_i-\dfrac{q+q'}{pk_0},
$$
i.e.,
$$
\alpha'_0=-\alpha_0+\sum_{i=1}^mt_i\alpha_i+N_2,
\qquad t_i,\,N_2 \in \mathbb{Z},\,\,\,i=1,\,2,\,\cdots,\,m.
$$

As a result, according to the above two cases and the known conditions
$$
\alpha_0=d\theta_0\,\,\,({\rm mod}\,\,\,\mathbb{Z}),
\qquad \alpha'_0=d\theta'_0\,\,\,({\rm mod}\,\,\,\mathbb{Z}),
\qquad \alpha_i=d\beta_i\,\,\,({\rm mod}\,\,\,\mathbb{Z}),
$$
where $i=1,\,2,\,\cdots,\,m$,
we have
$$
\theta'_0=\pm\,\theta_0+\sum_{i=1}^m\hbar_i\beta_i+\dfrac{N}{d},
\qquad \hbar_i,\,N \in \mathbb{Z},\,\,\,i=1,\,2,\,\cdots,\,m.
$$
Assume that $N \equiv N'\,\,\,({\rm mod}\,\,\,d)$.
Since ${\rm gcd}\,(d, k)=1$, then there exists some $\hbar_{m+1} \in \mathbb{Z}$ such that $\hbar_{m+1}k = N'\,\,\,({\rm mod}\,\,\,d)$.
So
$$
\theta'_0=\pm\,\theta_0+\sum_{i=1}^{m}{\hbar_i\beta_i}+\dfrac{N}{d}-\hbar_{m+1}\beta_{m+1}+\hbar_{m+1}\beta_{m+1}
=\pm\,\theta_0+\sum_{i=1}^{m+1}{\hbar_i\beta_i}\,\,\,({\rm mod}\,\,\,\mathbb{Z}),
$$
where $\beta_{m+1}=\dfrac{k}{d}$.
Take
$$
\theta_i=\beta_i,\,\,\,i=1,\,2,\,\cdots,\,m+1,
\qquad\theta_{i'}=0,\,\,\,i'=m+2,\,m+3,\,\cdots,\,n,
$$
and
$$
\ell_j=\hbar_j,\,\,\,j=1,\,2,\,\cdots,\,m+1,
\qquad \ell_{j'} \in \mathbb{Z},\,\,\,j'=m+2,\,m+3,\,\cdots,\,n.
$$
Therefore, we obtain
$$
\theta'_0=\pm\,\theta_0+\sum_{i=1}^n{\ell_i\theta_i}\quad({\rm{mod}}\,\,\,\mathbb{Z}),
\qquad \ell_i \in \mathbb{Z},\,\,\, i=1,\,2,\,\cdots,\,n.
$$

$\bm{(4)}$ Assume that $\rho(L_g)$ satisfies the conditions in {\bf Case 3},
and $\rho(L_{g'})$ satisfies the conditions in {\bf Case 1}.

Similar to the discussion in Case $(3)$, take left actions $L,\,L' \in \mathscr{M}_{T_{{\rm SU}(2) \times T^n}}$ satisfying
\begin{gather*}
\rho(L)=(\begin{array}{cccccccc}
\alpha_0&\alpha_1&\cdots&\alpha_m&0&\cdots&0
\end{array})^{\rm T},
\\
\rho(L')=(\begin{array}{cccccccc}
\alpha'_0&\alpha_1&\cdots&\alpha_m&0&\cdots&0
\end{array})^{\rm T},
\end{gather*}
where
$$
\alpha_0=d\theta_0\,\,\,({\rm mod}\,\,\mathbb{Z}),
\,\, \alpha'_0=d\theta'_0\,\,\,({\rm mod}\,\,\mathbb{Z}),
\,\,\alpha_i=d\beta_i\,\,\,({\rm mod}\,\,\mathbb{Z}),
\,\, i=1,\,2,\,\cdots,\,m.
$$
Then $L$ and $L'$ are topologically conjugate.
Let $h'$ be a topological conjugacy from $L$ to $L'$ satisfying $h'(e)=e$.
According to the known conditions, it is easy to see that $\alpha'_0$ is an rational number,
$\alpha_1,\,\alpha_2,\,\cdots\,\alpha_m$ are rational independent,
and $\alpha_0,\,\alpha_1,\,\cdots\,\alpha_m$ are $m+1$ rational dependent irrational numbers.
Assume that
\begin{align}\label{eq:7}
\left\{
\begin{array}{l}
k_0\alpha_0+k_1\alpha_1+\cdots+k_m\alpha_m+\dfrac{q}{p}=0,
\\[0.3 cm]
k'_0\alpha'_0+k'_1\alpha_1+\cdots+k'_m\alpha_m+\dfrac{q'}{p'}=0.
\end{array}
\right.
\end{align}
where $k_0$ is a positive integer, $k_0'=1$, $k_1'=k_2'=\cdots=k_m'=0$,
$p,\,p',\,q,\,q',\,k_i \in \mathbb{Z},\,i=1,\,2,\cdots,\,m$, and
$$
{\rm gcd}\,(k_0, k_1, \cdots, k_m)={\rm gcd}\,(p, q)={\rm gcd}\,(p', q')=1.
$$
Then by the same way in Case $\bm{(3)}$, we can obtain
$$
\theta'_0=\pm\,\theta_0+\sum_{i=1}^n{\ell_i\theta_i}\quad({\rm{mod}}\,\,\,\mathbb{Z}),
\qquad \ell_i \in \mathbb{Z},\,\,\, i=1,\,2,\,\cdots,\,n.
$$

In conclusion, the necessity of Theorem \ref{the:1} is proved

\end{proof}

At the end of this section, we give the proof the necessity of Proposition \ref{prop:6}.

\begin{proof}
Let $L_g$ and $L_{g'}$ be two left actions on ${\rm SU}(2) \times T^n$ and they are topologically conjugate.
According to the discussion in introduction, we know that there exist some $s,\,s' \in {\rm SU}(2) \times T^n$,
$t,\,t' \in T_{{\rm SU}(2) \times T^n}$ such that
$$
L_s \circ L_g=L_t \circ L_s,
\qquad L_{s'} \circ L_{g'}=L_{t'} \circ L_{s'},
$$
i.e. $L_g$ and $L_t$ are topologically conjugate, $L_{g'}$ and $L_{t'}$ are topologically conjugate,
and $L_s,\,L_{s'}$ are the topological conjugacies between them, respectively.
Thus, $L_t$ and $L_{t'}$ are topologically conjugate.
Assume that
$$
\rho(L_t)
=\left(\begin{array}{c}
\theta_0\\
\theta_1\\
\vdots\\
\theta_n
\end{array}\right),\qquad \rho(L_{t'})=\left(\begin{array}{c}
\theta'_0\\
\theta'_1\\
\vdots\\
\theta'_n
\end{array}\right),
$$
where $\theta_i,\,\theta'_i \in [0, 1),\,i=0,\,1,\,\cdots,\,n$.
Then it follows from the necessity of Theorem \ref{the:1} that we have
$$
\left(\begin{array}{c}
\theta'_0\\
\theta'_1\\
\vdots\\
\theta'_n
\end{array}\right)=\left(\begin{array}{cc}
\pm\, 1& \bm{\ell}\\
\bm{0}&{\bm A}
\end{array}\right)\left(\begin{array}{c}
\theta_0\\
\theta_1\\
\vdots\\
\theta_n
\end{array}\right)\quad(\rm{mod}\,\,\,\mathbb{Z}^{n+1}),
$$
where $\bm{\ell}=(\begin{array}{cccc}\ell_1&\ell_2&\cdots&\ell_n\end{array})$ is a $1 \times n$ integer matrix,
and ${\bm A} \in {{\rm GL}}_n(\mathbb{Z})$.
By the same way in $\S\,3.2$, we can construct the topological conjugacies from $L_t$ to $L_{t'}$,
and it is not difficult to see that these topological conjugacies are smooth homeomorphisms.
Hence $L_t$ and $L_{t'}$ are smooth conjugate.
Since the left actions on ${\rm SU}(2) \times T^n$ are all smooth homeomorphisms,
i.e., $L_s$ and $L_{s'}$ are smooth homeomorphisms,
so $L_g$ and $L_t$ are smooth conjugate,  $L_{g'}$ and $L_{t'}$ are smooth conjugate,
and then $L_g$ and $L_{g'}$ are smooth conjugate.
Therefore, the necessity of Proposition \ref{prop:6} is proved.

\end{proof}

\end{document}